# A 3D DPG MAXWELL APPROACH TO NONLINEAR RAMAN GAIN IN FIBER LASER AMPLIFIERS

S. NAGARAJ, J. GROSEK, S. PETRIDES, L. DEMKOWICZ, J. MORA

ABSTRACT. We propose a three dimensional Discontinuous Petrov-Galerkin Maxwell approach for modeling Raman gain in fiber laser amplifiers. In contrast with popular beam propagation models, we are interested in a truly full vectorial approach. We apply the ultraweak DPG formulation, which is known to carry desirable properties for high-frequency wave propagation problems, to the coupled Maxwell signal/pump system and use a nonlinear iterative scheme to account for the Raman gain. This paper also introduces a novel and practical full-vectorial formulation of the electric polarization term for Raman gain that emphasizes the fact that the computer modeler is only given a measured bulk Raman gain coefficient. Our results provide promising qualitative corroboration of the model and methodology used.

## 1. INTRODUCTION

The main aim of this paper is to present a full Maxwell, three dimensional (3D) Discontinuous Petrov-Galerkin (DPG) simulation of a fiber amplifier, using Raman gain [63, 42, 59, 44] in a typical passive, step-index, core-pumped optical fiber amplifier as the test case for initial validation purposes. In this regard, we present several novel advances, both in the modeling as well as in the methodology used to study Raman amplification in a full vectorial model.

First, our propagation model makes minimal assumptions on the electromagnetic fields in question, unlike the scalar beam propagation method (BPM, see [53, 58, 48, 64, 4] and references therein), which assumes a polarization maintaining propagation of the electromagnetic fields in an optical fiber, whereas our treatment is truly vectorial. Though both semi-vectorial and full vectorial BPM approaches have already been implemented (see [36, 37, 56, 57, 26] and references therein), we are introducing a fiber model that is a full boundary value problem rather than an initial value problem. In addition, we employ 3D isoparametric curvilinear elements to model the curved fiber (core and inner cladding) geometry, which can also later be used for studying microstructure fibers or hollow-core gas-filled fiber lasers. Indeed, most scalar fiber modeling techniques assume, starting with the initial condition, that only one of the three electric (and corresponding magnetic) field components dominate in magnitude during propagation, and thus treats the non-dominant components as zero. This is due to the assumption that the source light is robustly linearly polarized, and is thus only launched into one of the three electric field components, usually also neglecting the corresponding magnetic field component. Also, by assuming that the fiber is polarization maintaining (either by design or by active control), results in negligible field coupling as the light propagates through the fiber. While this assumption reduces the complexity of the model from a vectorial curl-curl Maxwell system to a scalar Helmholtz system, it may be the case that such assumptions may not hold to the degree required for the model to be accurate, especially in the presence of injected light that is not perfectly linearly polarized, or when there are high intensities, manufacturing defects, fiber bending, thermal effects, and/or the presence of embedded microstructures. In other words, the weakly coupled polarization states assumption may not be true in general, which would result in non-trivial coupling between the electromagnetic field components as the light propagates down the fiber.

Second, we propose a novel full vectorial time-harmonic 3D model for Raman gain. This paper shows how Raman gain, which typically is viewed as a nonlinear third-order susceptibility component of the electric polarization, can be derived by assuming that it originates from a mostly imaginary perturbation to the refractive index, just as active gain is usually derived. The proposed Raman model is particularly significant, since this fits well with, if not instrumental for, our full Maxwell simulation efforts, even though this effort centers on the validation of the numerical approach and not on a demonstration of polarization coupling.

Finally, our use of the DPG method, established by Demkowicz and Gopalakrishnan [18], for discretizing the fiber amplifier model equations is motivated by a host of desirable properties of the DPG method that have been theoretically and numerically explored, and validated in the recent past. The theoretical foundations of the subject were established in [16, 17, 10, 33, 11]. This active area of research has been successfully employed to problems in linear elasticity [38, 29], time harmonic wave propagation, including DPG versions of perfectly matched layers (PMLs) [54, 62], compressible and incompressible Navier-Stokes [25, 23, 12, 24], fluid flow [39], viscoelasticity [27] and space-time formulations [20, 35, 25, 22]. Moreover, versions of DPG for polygonal meshes have been introduced in [61]. Theoretical advances in goal-oriented adaptivity using DPG have been done in [41]. Practical implementation issues regarding conditioning of DPG systems are addressed in [40].

Indeed, the ideal DPG method (with optimal broken test functions) has been shown to provide a uniform, mesh-independent stable discretization for *any* well-posed variational formulation [16, 11]. The computationally tractable, practical DPG method [33], upon discretization of the so-called trial-to-test operator, retains the guaranteed stability with a numerically estimable stability constant [52]. The DPG method uses element-wise defined test spaces with no global conformity ("broken" test spaces), which allow for parallelism. Since the method can be recast as minimum residual, and also a mixed method with a built-in error indicator (the residual), one can have automatic *hp* adaptivity starting from an arbitrarily coarse mesh, which has importance in problems involving singularities. Finally, the method always delivers a sparse, Hermitian (symmetric), positive definite system making iterative conjugate gradient based solvers ideal for large systems that cannot be handled by direct solvers [54, 34, 5, 55].

Our model incorporates the fact that the ultraweak (UW) DPG formulation, used for solving the electromagnetic equations, provide us with both electric and magnetic fields. Thus, we are able to compute the cross sectional power via the time-averaged Poynting vector (irradiance) using the DPG trace variables. In this context, we also note the fact that we utilize a frequency domain PML (see [6, 60, 13, 31, 45, 46, 62]), which is also implemented using the ultraweak DPG formulation. As we shall explain, the use of a PML is critically important, and one cannot adequately observe the gain phenomenon with simpler impedance boundary conditions. Moreover, for element computations, we employ *sum factorization* to integrate the local DPG matrices, which significantly accelerates the otherwise temporally expensive element integration [47, 43].

Thus, to our knowledge, this paper is the first attempt at a general, full vectorial simulation of 3D Maxwell equations with a nonlinear gain term, equipped with a PML, in the context of higher-order Galerkin-based simulations. Also, this paper introduces an innovative formulation of Raman gain [63, 42, 59, 44] amenable to a vectorial simulation, which presumes that only the measured bulk Raman gain coefficient is available to the computer modeling team, as is almost always the case. Our simulations are, at this stage, not scaled to perform on supercomputing infrastructures, as would be needed for modeling fiber amplifiers of realistic sizes. However, the novelty of applying such a generalized approach to solving a vectorial, nonlinear fiber amplifier model with advanced 3D DPG technologies that provide the necessary accuracy and the unprecedented computational efficiency (for this type of methodology) is the major contribution of this paper to the literature. Subsequent investigations into code optimization, scaling and parallelism will allow for large-scale simulations of not only Raman gain but other phenomenon endemic to high-power fiber amplifiers such as stimulated Brillouin scattering (SBS), the transverse mode instability (TMI), thermal lensing, fiber bending, et cetera [2, 1, 7, 49, 50]. Also, this formulation of the governing equations, along with the resulting simulation approach, can serve as a basis for studying new amplifier configurations and/or for optimizing microstructure designs in future efforts.

We emphasize that our aim is to obtain *qualitative* results that indicate the feasibility of the methodology applied to the full Maxwell model, and, as such, we use an artificially large Raman gain coefficient so as to be able to see the gain effects in a short enough fiber length that all calculations can be accomplished on a single compute node. The remainder of this paper is organized as follows. In Section 2, we provide details of the physics underlying our model, introduce the novel, full vectorial electric polarization term that accounts for Raman gain, and delineate the system of equations that we will be solving. Section 3 briefly outlines how the DPG methodology is applied to general broken variational formulations. The discussions in Sections 2 and 3 are unified in Section 4, which provides details of the variational formulations, the nonlinear iterative scheme, and the time-harmonic Poynting theorem that are used in our model. We discuss in detail our results in Section 5, and conclude in Section 6. Three appendices provide details of the PML, the sum factorization implementation, and the theoretical underpinnings of the DPG approach.

**Acknowledgements**: This work was supported by AFOSR FA9550-17-1-0090 DOD Air Force Research Lab. We thank them for their support. Information in this paper is approved for public release on 08 May 2018 by AFRL OPSEC/PA OPS-18-19547. The work presented in this paper is drawn from the doctoral dissertation [51]. This work has been submitted to the Journal of Computational Physics (JCP).

## 2. 3D MAXWELL RAMAN GAIN MODEL

2.1. **Fiber Model.** This model considers a continuous wave (cw), double clad, non-dispersive, circularly symmetric, weakly guided, step-index fiber amplifier, where the core and cladding regions are isotropic and homogeneous (see Fig. 1). The outer layer (second cladding) of the fiber is a polymer coating that covers the inner cladding of fused silica. The refractive index of the core ($n_{core}$) and cladding ($n_{cladding}$) satisfy $n_{core} - n_{cladding} \ll 1$. Since it is assumed that all of the light (pump and signal laser fields) in this fiber is guided in the core region by total internal reflection, the subsequent model will ignore the polymer jacket, given that it has almost no effect on the core guided light.

Because both the pump and signal fields are seeded into the core region at the beginning of the fiber ($z = 0$), it is assumed that the light only propagates in the forward direction, which is a typical approximation for a co-pumped passive fiber amplifier. Such a configuration also suggests that both the pump and the signal are already highly coherent, which means that this amplifier acts only as a frequency converter, instead of also as a brightness enhancer.

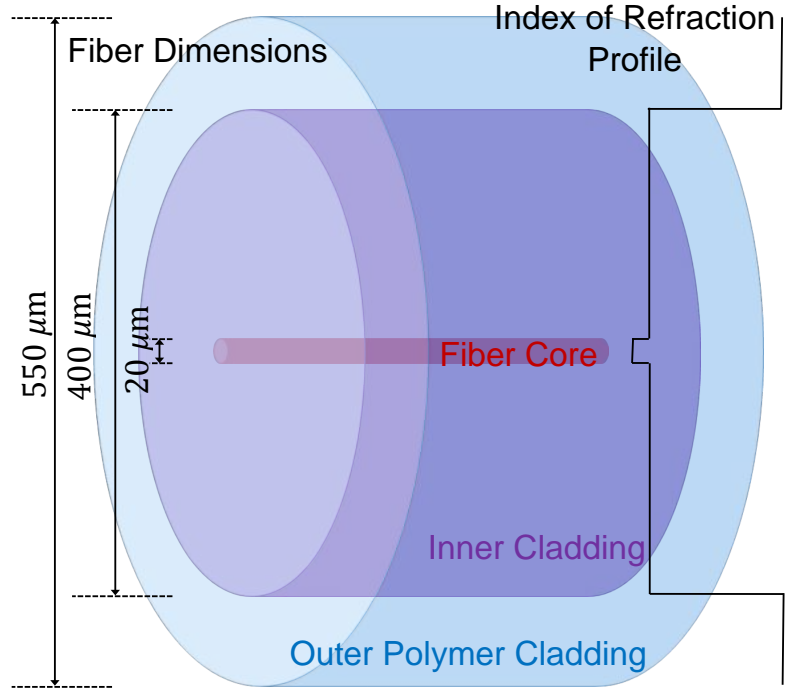

FIGURE 1. A typical circularly symmetric, double-clad, step index fiber amplifier, with a core region made of silica glass, a cladding region also made of silica glass, but with a slightly lower index of refraction than the core region, and a polymer coating, the outer cladding, with a substantially lower index of refraction than the inner cladding region. Such fibers are usually ∼5-100s meters long, even though it is depicted here as only being a few hundred microns long.

For the purposes of this Raman gain analysis, the electromagnetic fields are treated as time-harmonic. This is justified by the fact that real core-pumped Raman amplifiers are indeed usually seeded by lasers that produce near monochromatic light, and because any other sources of time dependent behaviour, most notably thermal effects, in passive fibers only occur at significantly slower varying time scales than the optical frequencies of the light present in the fiber. Thus, the following time-harmonic ansatz is assumed for all involved electromagnetic fields:

$$\mathbb{E}_0(x, y, z, t) = \mathbb{E}(x, y, z)e^{i\omega t} + \text{ c.c. and}$$
$$\mathbb{H}_0(x, y, z, t) = \mathbb{H}(x, y, z)e^{i\omega t} + \text{ c.c.},$$

where $\omega$ is the frequency of propagation, $i = \sqrt{-1}$ and c.c. indicates complex conjugate of the previous term.

The fact that the pump (p) and signal (s) fields are monochromatic and well-separated from one another, allows for solving two separate sets of Maxwell equations, which are coupled together through the Raman gain:

$$\nabla \times \mathbb{E}_l = -i\omega_l\, \mu \mathbb{H}_l,$$
$$\nabla \times \mathbb{H}_l = i\omega_l\, \varepsilon \mathbb{E}_l + i\omega_l \mathbb{P}_l,$$
$$\nabla \cdot \varepsilon \mathbb{E}_l = \rho,$$
$$\nabla \cdot \mu \mathbb{H}_l = 0,$$

(2.1)

where $l = \mathrm{p,s}$ is an index for the two frequencies of light, and $\mathbb{E}_l$ and $\mathbb{H}_l$ are the time-harmonic electric and magnetic fields respectively. The electric permittivity and magnetic permeability are denoted by $\varepsilon$ and $\mu$, respectively and the corresponding free-space quantities will be denoted by $\varepsilon_0$ and $\mu_0$. The electric charge density $\rho$ is zero for silica fibers, and $\mathbb{P}_l$ represents the electric polarization term.

### 2.2. Polarization Model.

Since silica fibers have negligible magnetic susceptibilities, all of the interactions between the electromagnetic fields and the medium can be formulated mathematically through the electric polarization term ($\mathbb{P}_l$). The electric polarization can be expanded in terms of the electric field and susceptibility tensors $\boldsymbol{\chi}^{(i)}, i = 1, 2, \ldots$ as follows...

$$\mathbb{P} = \varepsilon_0 \left( \underbrace{\boldsymbol{\chi}^{(1)} \cdot \mathbb{E}}_{\substack{\text{background refractive index (real)}\\ \text{active laser gain (imaginary)}}} + \boldsymbol{\chi}^{(2)} : \mathbb{E} \otimes \mathbb{E} + \underbrace{\boldsymbol{\chi}^{(3)} \vdots \mathbb{E} \otimes \mathbb{E} \otimes \mathbb{E}}_{\text{Raman gain} \propto |\mathbb{E}|^2 \mathbb{E}} + \ldots \right) [2].$$

An adequate model for this demonstration of a typical co-pumped passive fiber amplifier that experiences significant Raman gain must include the background index of refraction of the fiber, which will be denoted as $\mathbb{P}_l^{\text{background}}$ and is expressed through the real part of the first-order susceptibility. Also, the model must include the contribution of the Raman gain to the electric polarization, which will be denoted as $\mathbb{P}_l^{\text{Raman}}$ and is considered to be a component of the third-order susceptibility tensor. Active laser gain ($\mathbb{P}_l^{\text{active gain}}$) in a fiber amplifier is often seen as mostly imaginary perturbation to the refractive index, and is thus expressed as part of the first-order susceptibility term. This perturbation to the refractive index can be expressed as

$$\mathbf{n}_l^2 + 2\delta n^{\text{gain}} \frac{\mathbf{n}_l}{|\mathbf{n}_l|} \approx \frac{\varepsilon^l}{\varepsilon_0},$$

where $\boldsymbol{\varepsilon}^l$ is the dielectric tensor of the medium and $\delta n^{\text{gain}} = \delta n^{\text{gain}}(\omega_l)$ is a complex perturbation to the refractive index that causes a gain in the optical field. As will be shown presently, Raman gain can also be derived from the perspective that it is a mostly imaginary perturbation to the refractive index. A more complete model might include other effects such as linear loss ($\mathbb{P}_l^{\text{loss}}$), thermal effects ($\mathbb{P}_l^{\text{thermal}}$) and/or other optical nonlinearities ($\mathbb{P}_l^{\text{opt. nonlin.}}$) such as SBS, the Kerr nonlinearity, and/or four-wave mixing.

For the purposes of this paper, the electric polarization model takes the form of

$$\mathbb{P}_l(\mathbb{E}_l) = \mathbb{P}_l^{\text{background}}(\mathbb{E}_l) + \mathbb{P}_l^{\text{Raman}}(\mathbb{E}_l),$$

where

$$\mathbb{P}_l^{\text{background}}(\mathbb{E}_l) \approx \varepsilon_0 \big( \mathbf{n}_l^2 - \mathbb{I} \big) \mathbb{E}_l,$$

(2.2)

given that $\mathbb{I}$ is the identity tensor and $\mathbf{n}_l$ is the real-valued index of refraction tensor that accounts for the differences between the refractive indices of the fiber core region, the inner cladding region, and the polymer jacket region of the fiber [2, 7].

Raman scattering is an inelastic optical nonlinearity that occurs as incident light (the pump), at a sufficiently high-intensity, vibrates the molecules of the medium, resulting in optical phonons and scattered photons (the Stokes field), usually of a lower frequency than the incident photons. This process can start from noise, but in this model the Raman scattering is stimulated by having a seeded signal field offset in frequency from the pump field so as to achieve peak Raman gain and coinciding perfectly with the Stokes field frequency.

In order to derive the Raman gain contribution to the electric polarization, first consider how one might derive the contribution of active laser gain to the electric polarization. This approach is outlined in [63] using a scalar electric field; however, the process can be extended to a vectorial field. Even in high gain amplifiers, the gain is still a perturbation to the refractive index, and thus one should not expect that the gain would significantly contribute to the divergence of the electric field: $\boldsymbol{\nabla} \cdot \mathbb{P}_l^{\text{gain}} \approx 0$. The gain contribution to the first-order susceptibility can be denoted as $\boldsymbol{\chi}_{\text{g}}^{(1)} = \boldsymbol{\chi}_{\text{g}}^{(1)}(x, y, z, t)$ and can be decomposed into its real and imaginary components: $\boldsymbol{\chi}_{\text{g}}^{(1)}(x, y, z, t) = \boldsymbol{\chi}_{\text{g}}^{\text{Re}}(x, y, z, t) + i\boldsymbol{\chi}_{\text{g}}^{\text{Im}}(x, y, z, t)$, where $\boldsymbol{\chi}_{\text{g}}^{\text{Re/Im}}(x, y, z, t) \in \mathbb{R}$. It is reasonable to assume that the electric field grows according to a given gain function: $g_l = g_l(x, y, z, t)$, with units of $\text{m}^{-1}$. The electric field vector with gain is expressed as

$$\mathbb{E}_0^l(x, y, z, t) \approx \frac{1}{2}\mathbb{E}_l(x, y, z, t)e^{\frac{\langle g_l\rangle z}{2} + i(\omega_l t - \boldsymbol{\beta}^l \cdot \mathbf{r})} + \text{ c.c.},$$

where each component of propagation constant vector $\boldsymbol{\beta}^l$ is a positive real value, $\mathbf{r} = [x \ y \ z]^{\text{T}}$, the electric field envelop $\mathbb{E}_l$ is also slowly varying in time, $\omega_l > 0$, and

$$\langle g_l\rangle(z, t) = \frac{1}{\text{A}_{D_{\text{gain}}}} \iint_{D_{\text{gain}}} g_l(x, y, z, t) \ dxdy \text{ (by the 2D Mean Value Theorem)},$$

where $\text{A}_{D_{\text{gain}}}$ represents the transverse area of the domain of the gain. The expression for the electric field assumes that both the electric field amplitude ($\mathbb{E}_l$) and the gain function ($g_l$) are slowly varying compared to longitudinal oscillations at a frequency of $\beta_z^l$ and to temporal oscillations at a frequency of $\omega_l$. It is reasonably assumed that if the gain function obeys slowly varying envelope approximations than so does the gain contribution to the first-order susceptibility ($\boldsymbol{\chi}_{\text{g}}^{(1)}$). Therefore,

$$\begin{aligned}|\partial_{zz}f| \ll \beta|\partial_z f| \ll \beta^2|f| \\ |\partial_{tt}f| \ll \omega|\partial_t f| \ll \omega^2|f|\end{aligned} \text{ where } f \in \{\boldsymbol{\chi}_{\text{g}}^{\text{Re/Im}}, g, \mathbb{E}\}. \tag{2.3}$$

Furthermore, $\mathbb{E}_l$ is assumed to be slowly varying in $x$ and $y$ compared to the oscillations at the frequencies $\beta_x^l$ and $\beta_y^l$ in the $x$- and $y$-directions respectively.

It will be assumed that the vectorial Helmholtz equation robustly holds when applied to the slowly varying electric field amplitude:

$$\left[\boldsymbol{\Delta} + \frac{\mathbf{n}_l^2\omega_l^2}{c^2}\right]\left(\mathbb{E}_0^l e^{-\frac{\langle g_l\rangle z}{2}}\right) = 0,$$

basically indicating that the light propagates in the fiber even if there is no gain; i.e., the fiber is a waveguide. Moreover, it has been assumed that the light propagates only in the $z$-direction in order to simplify the mathematics, which means that Poynting vector $\mathbb{S}$, which parallel to the propagation constant vector $\boldsymbol{\beta} = \mathbf{n}_{\text{eff}}\omega/c$, is assumed to only have a $z$-component (and zero-valued $x$- and $y$-components). The gain has an isotropic effect, and thus its perturbative contribution to the refractive index occurs such that each direction of the refractive index is altered the same as any other direction even if the refractive index is birefringent (anisotropic). Mathematically, this is captured by $\boldsymbol{\chi}_{\mathrm{g}}^{\mathrm{Re/Im}} = \chi_{\mathrm{g}}^{\mathrm{Re/Im}} \mathbf{n}_l / |\mathbf{n}_l|$. Finally, the main idea of this derivation is to express the gain contribution to the electric polarization as a function of the first-order susceptibility due to gain:

$$\mathbb{P}_l^{\mathrm{gain}}(x,y,z,t) = \varepsilon_0 \left( \boldsymbol{\chi}_{\mathrm{g}}^{(1)}(x,y,z,t) \cdot \mathbb{E}_0^l(x,y,z,t) \right) = \varepsilon_0 \left( \boldsymbol{\chi}_{\mathrm{g}}^{\mathrm{Re}} + i\boldsymbol{\chi}_{\mathrm{g}}^{\mathrm{Im}} \right) \cdot \mathbb{E}_0^l. \tag{2.4}$$

Starting with the electric field wave equation in a dielectric medium with the gain contribution to the electric polarization term, and then applying the slowly varying envelope approximations and the vectorial Helmholtz equation:

$$\frac{\mathbf{n}_l}{|\mathbf{n}_l|} \left[ \boldsymbol{\Delta} \mathbb{E}_0^l - \frac{\mathbf{n}_l^2}{c^2} \frac{\partial^2 \mathbb{E}_0^l}{\partial t^2} \right] \approx \mu_0 \frac{\partial^2 \mathbb{P}_l^{\mathrm{gain}}}{\partial t^2} = \varepsilon_0 \mu_0 \frac{\partial^2}{\partial t^2} \left[ \left( \boldsymbol{\chi}_{\mathrm{g}}^{\mathrm{Re}} + i\boldsymbol{\chi}_{\mathrm{g}}^{\mathrm{Im}} \right) \cdot \mathbb{E}_0^l \right],$$

one derives,

$$\underbrace{\left( \frac{g_l}{2} \right)^2 + \left( \frac{\omega_l}{c} \right)^2 \chi_{\mathrm{g}}^{\mathrm{Re}}}_{\text{real-valued}} \approx \underbrace{i \left[ \left( \frac{\omega_l}{c} \right)^2 \chi_{\mathrm{g}}^{\mathrm{Im}} + g_l \beta_z^l \right]}_{\text{purely imaginary}}.$$

Now note that a real-valued function can only equal an imaginary-valued function when both functions are identically zero. Therefore, setting each side of the relation equal to zero produces the following relations:

$$\chi_{\mathrm{g}}^{\mathrm{Re}}(x,y,z,t) \approx - \left( \frac{g_l(x,y,z,t)c}{2\omega_l} \right)^2 \ \leftarrow \ \frac{g_l(x,y,z,t)c}{\omega_l} \ll 1 < n_l^2 \tag{2.5}$$

$$\chi_{\mathrm{g}}^{\mathrm{Im}}(x,y,z,t) \approx g_l(x,y,z,t) \beta_z^l \left( \frac{c}{\omega_l} \right)^2 \approx \frac{n_{\mathrm{eff}}^l c}{\omega_l} g_l(x,y,z,t) \tag{2.6}$$

Continuing the derivation, one finds that

$$\delta n_l^{\mathrm{gain}}(x,y,z,t) \approx \frac{ic g_l(x,y,z,t)}{2\omega_l} \approx -\frac{\sigma_l(x,y,z,t)}{\omega_l \varepsilon_0},$$

where $\sigma_l$ is the dielectric conductivity, which is another way of viewing gain in a fiber amplifier. This approximation also indicates that any contribution to the real component of the refractive index (2.5) by the presence of gain in the fiber is negligible in comparison to the imaginary component contribution (2.6). In fact, for light in the 1-2 $\mu$m wavelength range, the perturbation to the real component of the index of refraction due to gain is about 7 orders of magnitude smaller than it is for the imaginary component. Therefore, one can approximate $\chi_{\mathrm{g}}^{\mathrm{Re}}(\mathbf{r}, t)$ to be zero, and re-express the contribution of gain to the electric polarization (2.4) using the derived relation for the imaginary component of the susceptibility due to gain (2.6) to get,

$$\mathbb{P}_l^{\mathrm{gain}}(x,y,z,t) \approx \frac{i\varepsilon_0 c \mathbf{n}_l}{\omega_l} g_l(x,y,z,t) \mathbb{E}_l(x,y,z,t). \tag{2.7}$$

Unfortunately, experimentalists consistently measure the bulk Raman gain coefficient as the primary means of determining how susceptible a fiber may be to experiencing the onset of Raman scattering. This means that simulations cannot produce a model better than the limitations imposed by this constant, and the methodology used to determine its value. It is important to understand that experimentalists ascertain the Raman gain coefficient measurement from a coupled set of ODEs for power ($P_l$) of the pump ($l = \mathrm{p}$) and Raman Stokes ($l = \mathrm{S}$) fields along the length of the fiber, which can be derived from multiple simplifying assumptions applied to Maxwell's equations [42, 59, 44]. Written concisely, without including extra terms for starting the Raman scattering from noise, this coupled set of ODEs takes the form of

$$\frac{dP_\mathrm{P}}{dz}(z) = \frac{\Upsilon_\mathrm{R}^\mathrm{p} g_\mathrm{R}}{\mathrm{A_{eff}}} P_\mathrm{P}(z) P_\mathrm{S}(z) \text{ and } \frac{dP_\mathrm{S}}{dz}(z) = \frac{\Upsilon_\mathrm{R}^\mathrm{S} g_\mathrm{R}}{\mathrm{A_{eff}}} P_\mathrm{P}(z) P_\mathrm{S}(z) \text{ with}$$
$$\mathrm{A_{eff}} = \frac{\iint_{\mathrm{A_{clad}}} (\varphi^\mathrm{p})^2 \, dxdy \iint_{\mathrm{A_{clad}}} (\varphi^\mathrm{S})^2 \, dxdy}{\iint_{\mathrm{A_{clad}}} (\varphi^\mathrm{p})^2 (\varphi^\mathrm{S})^2 \, dxdy}, \tag{2.8}$$

where $g_\mathrm{R} \in \mathbb{R}^+$ is the measured bulk Raman gain coefficient, and $\varphi^l = \varphi^l(x,y)$ is a single transverse mode of the fiber; presumably the fundamental mode. The dimensionless parameter $\Upsilon_\mathrm{R}^l$, with $l \in \{\mathrm{p}, \mathrm{S}\}$, allows for photon flux conservation when

$$\Upsilon_\mathrm{R}^l = \begin{cases} \frac{-\omega_\mathrm{p}}{\omega_\mathrm{S}}, & \text{when } l = \mathrm{p} \\ 1, & \text{when } l = \mathrm{S} \end{cases}.$$

These ODEs assume that the light, at either frequency, only resides in one transverse mode, which may be a limiting factor when considering large mode area (LMA) fibers. Also, note that the power at a particular point along the fiber is already independent of the transverse direction, and that the effective area calculation further washes out any transverse dependencies. Finally, recall that, in this simulation, the pump frequency is higher than the Stokes (or signal) frequency: $\omega_\mathrm{p} > \omega_\mathrm{S}$, or equivalently, $\lambda_\mathrm{S} > \lambda_\mathrm{p}$, so that the term $|\frac{-\omega_\mathrm{p}}{\omega_\mathrm{S}}| > 1$, which results in energy transfer from the pump field to the Stokes (signal) field.

A slight generalization to these power evolution ODEs (2.8) is the set of transverse dependent PDEs for the evolution of the irradiance ($I_l := |\mathrm{Real}(\mathbb{E}_l \times \mathbb{H}_l^*)|$) along the fiber:

$$\frac{\partial I_\mathrm{P}}{\partial z} = \Upsilon_\mathrm{R}^\mathrm{p} g_\mathrm{R} I_\mathrm{P} I_\mathrm{S} \text{ and } \frac{\partial I_\mathrm{S}}{\partial z} = \Upsilon_\mathrm{R}^\mathrm{S} g_\mathrm{R} I_\mathrm{P} I_\mathrm{S}. \tag{2.9}$$

These PDEs help illuminate the gain function ($g_l$) for Raman scattering.

Accepting that the bulk Raman gain coefficient is the primary means of determining Raman gain, it is prudent to introduce this constant directly into the derivation of gain; specifically by including $g_\mathrm{R}$ into the gain function ($g_l$). The gain function for Raman scattering can be extracted from the coupled irradiance PDEs (2.9) by choosing $g_l(x,y,z,t) = \Upsilon_\mathrm{R}^l g_\mathrm{R} I_k(x,y,z,t)$, where $k \neq l \in \{\mathrm{p}, \mathrm{s}\}$ and $\mathrm{S} \equiv \mathrm{s}$ for this simulation. Using this form of the Raman gain function in the expression for the contribution of gain to the electric polarization (2.7), which meets the necessary criteria of having units of $\mathrm{m}^{-1}$ and obeying the slowly varying envelop approximations (2.3), yields a novel and practical formulation of the Raman gain contribution to the electric polarization:

$$\mathbb{P}_l^\mathrm{Raman}(\mathbb{E}_l) \approx \frac{i\varepsilon_0 \mathbf{n}_l c}{\omega_l} (\Upsilon_\mathrm{R}^l g_\mathrm{R} I_k) \mathbb{E}_l. \tag{2.10}$$

Recalling that the intensity can be related to the square of the electric field (and likewise with the irradiance) it is clear that this expression sets $\mathbb{P}_\mathrm{Raman} \propto |\mathbb{E}|^2 \mathbb{E}$ as would be expected for a third-order susceptibility component of the electric polarization since $\mathbb{P}_\mathrm{Raman}^l \propto I_k \mathbb{E}_l$ and $I_k \propto |\mathbb{E}|_k^2$. These expressions show that the Raman gain is the source of the nonlinearity for this model of a coupled system of Maxwell equations.

2.3. **Non-dimensionalization of Governing Equations.** The above mentioned equations are dimensional. The non-dimensional version of the equations are derived in order to distinguish the physics from the system of units, especially since fibers have very disparate geometric scales. Indeed, the physical dimensions of a typical high-power fiber amplifier are 1 mm in diameter (at most) and about 5-100 m in length, but possibly even longer for some Raman amplifiers.

Let $l_0$ be a generic spatial scaling such that $x = l_0\hat{x}$, where $\hat{x}$ is the non-dimensional spatial variable. For the rest of the paper, the hat symbol ($\hat{\ }$) will indicate a non-dimensional parameter or variable. One can now derive that

$$\frac{\partial}{\partial x} = \frac{\partial}{\partial \hat{x}}\frac{\partial \hat{x}}{\partial x} = \frac{1}{l_0}\frac{\partial}{\partial \hat{x}}.$$

The corresponding non-dimensional curl operator will be denoted by $\hat{\nabla}\times$. Now consider the dimensionless versions of the electromagnetic fields, and frequency parameter, which can be expressed as

$$\mathbb{E}_l = E_0\hat{\mathbb{E}}_l, \ \mathbb{H}_l = H_0\hat{\mathbb{H}}_l, \text{ and } \omega_l = \omega_0\hat{\omega}_l.$$

In this simulation, the frequencies of the pump and signal fields are both near to $\omega_0 = 10^{15}$ rads/sec, which will be considered as the chosen value for this parameter. The other parameters that have been introduced for the non-dimensionalization of the governing equations will be chosen as follows:

$$l_0 = \frac{c}{\omega_0}, \ H_0 = \frac{1}{c}\sqrt{\frac{\omega_0\kappa_a}{\mu_0 g_R}}, \text{ and } E_0 = \sqrt{\frac{\mu_0\omega_0\kappa_a}{g_R}}.$$

With these choices, and with the identity $c^2 = (\varepsilon_0\mu_0)^{-1}$, one can determine that

$$\frac{\omega_0\mu_0 H_0 l_0}{E_0} = 1, \ \frac{\omega_0\varepsilon_0 E_0 l_0}{H_0} = 1, \text{ and } \frac{cg_R E_0 H_0}{\omega_0} = \kappa_a.$$

In order to augment the Raman gain phenomenon within a short fiber (several tens of wavelengths), we have introduced the artificial scaling parameter $\kappa_a$. This non-physical parameter scales the intensity values (through boundary conditions) throughout the fiber, and thereby allows us to simulate very short fiber lengths while allowing for the gain of signal power from the pump field. Now the first-order Maxwell system (2.1), with the expressions for the background refractive index (2.2) and Raman gain (2.10) contributions to the electric polarization, can be non-dimensionalized:

$$\begin{aligned} \hat{\nabla}\times\hat{\mathbb{E}}_l &= -i\hat{\omega}_l\hat{\mathbb{H}}_l \\ \hat{\nabla}\times\hat{\mathbb{H}}_l &= i\hat{\omega}_l\hat{\mathbb{E}}_l + i(\mathbf{n}_l^2 - \mathbb{I})\hat{\omega}_l\hat{\mathbb{E}}_l + i\hat{\omega}_l\frac{i\mathbf{n}_l\kappa_a\Upsilon_R^l}{\hat{\omega}_l}\big|\mathrm{Real}(\hat{\mathbb{E}}_k\times\hat{\mathbb{H}}_k^*)\big|\hat{\mathbb{E}}_l, \end{aligned} \tag{2.11}$$

which result in:

$$\begin{aligned} \hat{\nabla}\times\hat{\mathbb{E}}_l &= -i\hat{\omega}_l\hat{\mathbb{H}}_l \\ \hat{\nabla}\times\hat{\mathbb{H}}_l &= i\mathbf{n}_l^2\hat{\omega}_l\hat{\mathbb{E}}_l - \mathbf{n}_l\kappa_a\Upsilon_R^l\big|\mathrm{Real}(\hat{\mathbb{E}}_k\times\hat{\mathbb{H}}_k^*)\big|\hat{\mathbb{E}}_l \end{aligned} \tag{2.12}$$

This represents two nonlinear Maxwell systems, one with $l = \mathrm{s}$ and $k = \mathrm{p}$ and the other with $l = \mathrm{p}$ and $k = \mathrm{s}$, that are coupled together through the Raman gain term.

2.4. **Boundary conditions.** In this model, the light propagates along the $z$-axis in the fiber core, only in the forward ($+z$) direction. Since the model is formulated as a boundary value problem, it is paramount that the boundary conditions, especially on the output facet of the fiber, correctly capture the physics of the amplifier. The fiber is excited at the input end (corresponding to $z = 0$) with two light sources introduced into the fiber core region. Recall that we have introduced the artifical scaling parameter $\kappa_a$. The non-dimensionalization relations show that by choosing $0 < \kappa_a < 1$, we artifically increase the field intensities, thereby injecting an increased amount of power within the short fiber at $z = 0$ in order to see sufficient gain in a short distance along the fiber. A zero boundary condition for the electromagnetic fields is set at the outer edge of the inner cladding, which is far enough away from the core region so as to not significantly affect the guided light. Indeed, for guided light, the fields decay exponentially within the cladding, and at the radial boundary, $\sqrt{x^2 + y^2} = r = r_{\text{cladding}}$, the fields are, within numerical precision, zero. Finally, at the exit end of the fiber ($z = L$, where $L$ is the length of the fiber), appropriate out-flowing radiation boundary conditions are set. In order to facilitate this, a PML is introduced at the end of the fiber. The need for this is better understood by observing that the gain polarization can be thought of as producing an electric conductivity within the material.

2.4.1. *Implication of Conductivity on Boundary Condition:* How the Raman gain can be viewed in terms of non-zero conductivity $\sigma$ will now be addressed. We will drop the "hats" while referring to the non-dimensional equations derived earlier. Consider the term $i\omega_l \mathbb{P}_l$. The background part of the polarization behaves linearly:

$$i\omega_l \mathbb{P}_l^{\text{background}}, (\mathbb{E}_l) = i\omega_l(\mathbf{n}_l^2 - \mathbb{I})\mathbb{E}_l.$$

However, the Raman term yields:

$$i\omega_l \mathbb{P}_l^{\text{Raman}}(\mathbb{E}_l) = i\omega_l \frac{i\mathbf{n}_l}{\omega_l} \kappa_a \, \Upsilon_R^l |\text{Real}(\mathbb{E}_k \times \mathbb{H}_k^*)| \mathbb{E}_l = -(\mathbf{n}_l \kappa_a \Upsilon_R^l)|\text{Real}(\mathbb{E}_k \times \mathbb{H}_k^*)|\mathbb{E}_l.$$

The term $(\mathbf{n}_l \kappa_a \Upsilon_R^l)|\text{Real}(\mathbb{E}_k \times \mathbb{H}_k^*)|$ is purely real and hence can be interpreted as a material conductivity, which acts as a nonlinear coupling between the signal and pump fields. The entire amplification properties hinge on this term. Indeed, this nonlinear term is responsible for the power transfer from the pump field into the signal field, since $\Upsilon_R^s = 1$ while $\Upsilon_R^p = -\frac{\omega_p}{\omega_s} < -1$, which implies loss from the pump into the signal. The nonlinearity is, however, just a perturbation to the linear problem. Although the numerical value of the gain is increased significantly, it is a weak nonlinearity, since it does not induce any self-coupling in the signal and pump fields individually.

One implication for DPG implementation is apparent: a simple impedance-like boundary condition at the terminal end of the fiber will not suffice. Indeed, impedance boundary conditions for waveguides operate on the principle of a single propagating mode in a lossless linear medium with an exactly known impedance constant, say $\gamma$. One then relates the $\mathbb{E}, \mathbb{H}$ fields on a boundary (with normal $\vec{n}$) as:

$$\mathbb{E} + \gamma \, \vec{n} \times \mathbb{H} = 0.$$

However, since this is a nonlinear problem where the conductivity changes along the length of the fiber, an exact impedance-like relation between the $\mathbb{E}, \mathbb{H}$ on the terminal boundary is inapplicable and would lead to incorrect boundary behaviour resulting in a physically meaningless solution. Thus, one must develop a PML at the exit end of the fiber, which would not hamper the behaviour of the fields within the domain. Towards this end, PMLs have been widely used in finite element implementations. Most notably, [8] and the recent work [62] use DPG methods to implement ultraweak formulations for various wave propagation phenomena. In this case, a stretched coordinate PML for the ultraweak formulation is used, with stretching along the $z$-axis, since outgoing waves need to be attenuated in only the $z$-direction. We refer the interested reader to [62] and Appendix A for implementation details. Figure 2 indicates the use of ultraweak DPG PML for a rectangular waveguide with a pronounced exponentially growing wave. Notice how the wave attenuates completely after entering the PML region. A similar approach is used for the fiber geometry.

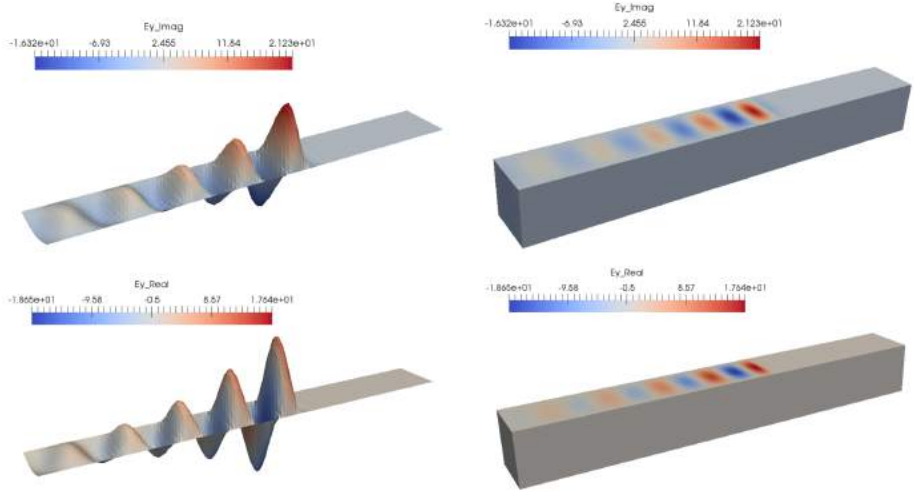

FIGURE 2. Ultraweak DPG PML for growing waves

Given the nature of this simulation, and the computational challenges it entails, the aim is to demonstrate *qualitative* results of the Raman gain action. This goal, for now, requires the use of a sufficiently short fiber so that all calculations can be completed on a single compute node in a reasonable amount of time. This is done with the understanding that future efforts will parallelize this model and implement it on a supercomputing platform, where more realistically sized fibers can be studied. Therefore, this simulation sets the fiber length to be less than 0.1 mm (~50-100 wavelengths), and artificially increases the field intensities (and thereby powers) in the fiber by many orders of magnitude in order to absorb significant amounts of the pump field in this short distance, allowing one to qualitatively observe the Raman process.

## 3. DPG Technology

The DPG technology is a multi-faceted approach to the stable discretization of well-posed variational formulations. In essence, DPG methods (used with optimal test functions) come with several impressive properties: uniform, mesh independent stability, localizable test norms via broken test spaces and a built-in canonical error indicator.

This section is devoted to providing a brief overview of the DPG methodology for any abstract variational formulation. This is motivated by the fact that the DPG method is applicable to any well-posed variational formulation, and its key properties are best described in an abstract setting. The second half of this section deals with energy spaces that arise in the discretization of the time-harmonic Maxwell system, and there the so-called ultraweak variational formulation, which shall be used for this Raman gain model, is defined. The following section will specialize the discussions of this section to the case of the nonlinear signal/pump time-harmonic Maxwell systems that arise in the context of the Raman gain problem.

In short, the DPG approach to discretization is as follows. First, identify a well-posed *broken* variational formulation corresponding to the physical problem of interest. Next, fix a finite dimensional subspace of the desired solution (trial) space, and compute the corresponding finite dimensional *optimal test space* by means of inverting the *enriched Riesz map*. Finally, the resulting (provably stable) discrete variational formulation is cast as a linear system and solved either by direct or iterative means. Refer to Appendix C for a more elaborate discussion of DPG theory.

### 3.1. **DPG in a Nutshell.**

Consider a variational formulation consisting of a continuous bilinear (or sesquilinear) form $\mathfrak{b}(\cdot, \cdot)$ defined on the product $\mathfrak{X} \times \mathfrak{Y}$ of Hilbert spaces $\mathfrak{X}, \mathfrak{Y}$ and a linear (or anti-linear) form $\mathfrak{l}(\cdot)$ defined on $\mathfrak{Y}$:

$$\text{find } \mathfrak{u} \in \mathfrak{X} \text{ such that } \mathfrak{b}(\mathfrak{u}, \mathfrak{v}) = \mathfrak{l}(\mathfrak{v}),$$

for all $\mathfrak{v} \in \mathfrak{Y}$. The symbols $\mathfrak{X}, \mathfrak{Y}$ are referred to as the trial and test spaces respectively. Such variational formulations arise naturally via the relaxation (i.e., integrating by parts) of the governing equations [15]. Only well-posed (in the sense of Hadamard [11, 9, 15]) variational formulations are of interest. Discretization of a well-posed variational formulation amounts to choosing finite dimensional subspaces $\mathfrak{X}_h, \mathfrak{Y}_h$ of $\mathfrak{X}, \mathfrak{Y}$ respectively and considering the discrete variational problem:

$$\text{find } \mathfrak{u}_h \in \mathfrak{X}_h \text{ such that } \mathfrak{b}(\mathfrak{u}_h, \mathfrak{v}_h) = \mathfrak{l}(\mathfrak{v}_h),$$

for all $\mathfrak{v} \in \mathfrak{Y}_h$. It can be shown (see [3]) that arbitrary choices of finite dimensional trial and test spaces do not yield well-posed discrete formulations. The *ideal DPG method* [18, 32] is a specialized discretization procedure which always guarantees a stable discretization by computing the so-called "optimal test functions" (see Appendix C and references therein). Distressingly, the optimality properties of the ideal DPG method can be computationally intractable in practice. Indeed, the exact computation of optimal test functions involves inverting the Riesz operator on the test space $\mathfrak{Y}$, which is an infinite-dimensional optimization problem. In practice, however, one truncates the inversion of the Riesz operator to a large, yet finite dimensional "enriched" test space $\mathfrak{Y}_r$. This truncation has, thankfully, a very benign effect on the overall stability. Construction of appropriate *Fortin operators* show the well-posedness and stability of the truncated ("practical") DPG computations [52, 33, 11]. Thus, practical implementations of DPG are optimal.

### 3.2. **Energy Spaces for the Maxwell Equations.**

Returning to the DPG discretization of time-harmonic Maxwell equations, there are, like other equations of physics [15], four conceivable variational formulations of the Maxwell equations [11]. In particular, two formulations stand out: the primal and ultraweak. As described in Appendix C, the ultraweak formulation is the formulation of choice in this Maxwell fiber amplifier problem. The remainder of this subsection defines the energy spaces required for the Maxwell system and defines the ultraweak formulation.

Consider a bounded, simply connected Lipschitz domain $\Omega \subset \mathbb{R}^3$ with boundary $\partial\Omega$ and unit normal vector $n$. The existence of a mesh $\Omega_h$ of finitely many open elements $K$, each with unit normal $n_K$, such that $\overline{\Omega} \subset \bigcup_{K \in \Omega_h} \overline{K}$ is assumed. Next, define:

$$
\begin{aligned}
L^2(\Omega) &:= \{f : \Omega \to \mathbb{R} : \int_\Omega |f|^2 < \infty\}, \\
\mathbb{L}^2(\Omega) &:= L^2(\Omega) \times L^2(\Omega) \times L^2(\Omega), \\
H(\mathrm{curl}, \Omega) &:= \{\mathbb{E} \in \mathbb{L}^2(\Omega) : \nabla \times \mathbb{E} \in \mathbb{L}^2(\Omega)\}, \\
H_0(\mathrm{curl}, \Omega) &:= \{\mathbb{E} \in H(\mathrm{curl}, \Omega) : n \times \mathbb{E}|_{\partial\Omega} = 0\}.
\end{aligned}
\tag{3.1}
$$

The broken counterpart of $H(\mathrm{curl}, \Omega)$ is defined as:

$$
H(\mathrm{curl}, \Omega_h) := \{\mathbb{E} \in \mathbb{L}^2(\Omega) : \mathbb{E}|_K \in H(\mathrm{curl}, K), K \in \Omega_h\} = \prod_{K \in \Omega_h} H(\mathrm{curl}, K). \tag{3.2}
$$

Notice that the broken counterpart of $\mathbb{L}^2(\Omega)$ is itself. The element-wise summed $\mathbb{L}^2(\Omega)$ inner product of the two arguments is denoted by $(\cdot, \cdot)_h$, and the element-wise summed duality pairing of appropriate dual spaces is represented by $\langle \cdot, \cdot \rangle_h$. The symbol $\| \cdot \|$ shall mean the $\mathbb{L}^2(\Omega)$ norm. As was shown in [11], the definition of trace operators are required in order to elegantly define the DPG interface spaces. First, define element trace operators:

$$
\begin{aligned}
\mathrm{t}_{K,\top}(\mathbb{E}) &:= (n_K \times \mathbb{E}) \times n_K|_{\partial K} \\
\mathrm{t}_{K,\perp}(\mathbb{E}) &:= (n_K \times \mathbb{E})|_{\partial K}
\end{aligned}
\tag{3.3}
$$

Notice that these trace operators have range $H^{-1/2}(\mathrm{curl}, \partial K)$ and $H^{-1/2}(\mathrm{div}, \partial K)$ respectively, i.e.,

$$
\begin{aligned}
\mathrm{t}_{K,\top} &: H(\mathrm{curl}, K) \to H^{-1/2}(\mathrm{curl}, \partial K), \\
\mathrm{t}_{K,\perp} &: H(\mathrm{curl}, K) \to H^{-1/2}(\mathrm{div}, \partial K).
\end{aligned}
\tag{3.4}
$$

Finally, the trace operators on the full broken $H(\mathrm{curl}, \Omega_h)$ space are defined via the element-wise application of the element trace operators:

$$
\begin{aligned}
\mathfrak{T}_\top &: H(\mathrm{curl}, \Omega_h) \to \prod_{K \in \Omega_h} H^{-1/2}(\mathrm{curl}, \partial K), \\
\mathfrak{T}_\perp &: H(\mathrm{curl}, \Omega_h) \to \prod_{K \in \Omega_h} H^{-1/2}(\mathrm{div}, \partial K).
\end{aligned}
\tag{3.5}
$$

The operators $\mathfrak{T}_\top, \mathfrak{T}_\perp$ are linear by construction. Finally, the spaces of interface variables (or interface spaces) can be defined as the images under the trace maps of the conforming $H(\mathrm{curl}, \Omega)$ space:

$$
\begin{aligned}
H^{-1/2}(\mathrm{div}, \partial\Omega_h) &:= \mathfrak{T}_\top(H(\mathrm{curl}, \Omega)) \\
H^{-1/2}(\mathrm{curl}, \partial\Omega_h) &:= \mathfrak{T}_\perp(H(\mathrm{curl}, \Omega)).
\end{aligned}
\tag{3.6}
$$

As shown in [11], the trace (quotient) norms on the two interface spaces are dual to each other.

### 3.3. Ultraweak Variational Formulation.
The ultraweak formulation (see also Appendix C) corresponds to the case where

$$
\mathfrak{X}_0 = \mathbb{L}^2(\Omega) \times \mathbb{L}^2(\Omega), \hat{\mathfrak{x}} = H^{-1/2}(\mathrm{curl}, \partial\Omega_h) \times H^{-1/2}(\mathrm{curl}, \partial\Omega_h),
$$

$$
\mathfrak{Y}_0 = H(\mathrm{curl}, \Omega) \times H_0(\mathrm{curl}, \Omega), \mathfrak{Y} = H(\mathrm{curl}, \Omega_h) \times H(\mathrm{curl}, \Omega_h).
$$

Denote by $\mathfrak{u} = (\mathbb{E}, \mathbb{H}) \in \mathfrak{X}_0$, $\hat{\mathfrak{u}} = (\hat{\mathbb{E}}, \hat{\mathbb{H}}) \in \hat{\mathfrak{x}}$ and $\mathfrak{v} = (\mathbb{R}, \mathbb{S}) \in \mathfrak{Y}$. The bilinear forms corresponding to the ultraweak formulation are:

$$
\begin{aligned}
\mathfrak{b}_0(\mathfrak{u}, \mathfrak{v}) &= (\mathbb{H}, \nabla \times \mathbb{R})_h - (i\omega\epsilon + \sigma)(\mathbb{E}, \mathbb{R})_h + (\mathbb{E}, \nabla \times \mathbb{S})_h + i\omega\mu(\mathbb{H}, \mathbb{S})_h, \\
&= (\mathfrak{u}, A^*\mathfrak{v})_h, \\
\hat{\mathfrak{b}}(\hat{\mathfrak{u}}, \mathfrak{v}) &= \langle n \times \hat{\mathbb{H}}, \mathbb{R} \rangle_h + \langle n \times \hat{\mathbb{E}}, \mathbb{S} \rangle_h.
\end{aligned}
\tag{3.7}
$$

Note that in the above expressions, we succinctly write the bilinear form $\mathfrak{b}_0(\cdot,\cdot)$ in terms of the (formal) adjoint $A^*$ of the Maxwell operator. The ultraweak formulation comes equipped with the (scaled) adjoint graph norm:

$$\|\mathfrak{v}\|^2_{\mathfrak{V}} := \alpha\|\mathfrak{v}\|^2 + \|A^*\mathfrak{v}\|^2.$$

Modification of the true adjoint graph norm (consisting of only $\|A^*\mathfrak{v}\|^2$ term) by adding the above $\alpha\|\mathfrak{v}\|^2$ scaling term is required to make the norm localizable [19, 17]. We use $\alpha = 1$. Next, the ultraweak DPG discretization of the Raman gain problem will be considered.

## 4. Setup of Simulations

Having established the superiority of the ultraweak formulation in Appendix C, this formulation will be used for the remainder of this work. For notational convenience, the "hats" that denoted the non-dimensional quantities derived in Section 2 shall be omitted.

### 4.1. Model Implementation.
For the Raman fiber amplifier simulations, shape functions developed in [28] are used, which support 3D elements of all shapes (hexahedron, prism, tetrahedron and pyramid). The coding for this problem was done on the $hp3D$ infrastructure detailed in the book [21]. As is noted in [40], the DPG method can be implemented in any standard finite element code. An all-hexahedron mesh is used in order to take advantage of the fast quadrature developed in [47]. The space of polynomials of order $p$ is denoted by $\mathcal{P}^p$, with $\mathcal{Q}^{(p,q,r)} := \mathcal{P}^p \otimes \mathcal{P}^q \otimes \mathcal{P}^r$ and

$$\begin{aligned}
\mathcal{W}_p &:= \mathcal{Q}^{(p,q,r)}, \\
\mathcal{Q}_p &:= \mathcal{Q}^{(p-1,q,r)} \times \mathcal{Q}^{(p,q-1,r)} \times \mathcal{Q}^{(p,q,r-1)}, \\
\mathcal{V}_p &:= \mathcal{Q}^{(p,q-1,r-1)} \times \mathcal{Q}^{(p-1,q,r-1)} \times \mathcal{Q}^{(p-1,q-1,r)}, \\
\mathcal{Y}_p &:= \mathcal{Q}^{(p-1,q-1,r-1)}.
\end{aligned} \tag{4.1}$$

The Maxwell system utilizes the Nedelec hexahedron of second type characterized by the exact sequence [21]:

$$\mathbb{R} \xrightarrow{\text{id}} \mathcal{W}_p \xrightarrow{\nabla} \mathcal{Q}_p \xrightarrow{\nabla\times} \mathcal{V}_p \xrightarrow{\nabla\cdot} \mathcal{Y}_p \longrightarrow 0 \ .$$

In the process of coding the Raman problem within $hp3D$, separate data structures for both signal and pump variables are supported, but the memory is allocated for each solve separately. This is possible due to the weak coupling between the two sets of fields through the Raman gain. Thus, while solving for the signal fields, memory is allocated only for the signal, and likewise while solving for the pump. The solvers used in this work come from the MUMPS (MUltifrontal Massively Parallel sparse direct Solver, see at http://mumps.enseeiht.fr/) library and the Intel MKL Pardiso solver.

### 4.2. Model Parameters.
This test problem of a core-pumped, step-index Raman amplifier sets the (non-dimensionalized) core and inner cladding radii to $\hat{r}_{\text{core}} = 0.25\sqrt{2}$ and $\hat{r}_{\text{cladding}} = 2.5\sqrt{2}$ respectively, and sets the cladding refractive index to $n_{\text{cladding}} = 1.45$. Using a numerical aperture of NA $\approx 0.0659$, and knowing that NA $= \sqrt{n_{\text{core}}^2 - n_{\text{cladding}}^2}$, the core refractive index can be determined as $n_{\text{core}} \approx 1.4515$. This means that the normalized frequency (or V-number) of the fiber is

$$\text{V} = \frac{2\pi r_{\text{core}}}{\lambda}\text{NA} \approx 2.198.$$

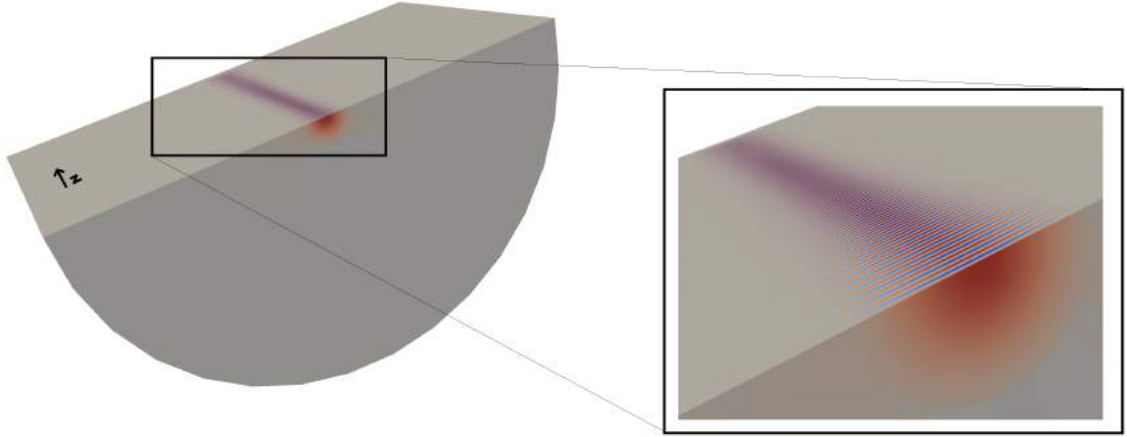

FIGURE 3. Cross sectional view of the solution to the linear problem with only background polarization for a fiber of length ≈ 80 wavelengths. The core is discretized with 5 hexahedral elements while the cladding has 4 hexahedral elements in the initial mesh. The zoomed part shows the close-up of the core region of the fiber. Here, the core radius is roughly one tenth the radius of the cladding.

Note that the V-number can also be expressed in terms of non-dimensional quantities as:

$$\mathrm{V} = \frac{\omega}{c} r_{\mathrm{core}} \mathrm{NA} = \frac{l_0 \omega_0}{c} \hat{\omega} \hat{r}_{\mathrm{core}} \mathrm{NA} = \hat{\omega} \hat{r}_{\mathrm{core}} \mathrm{NA},$$

where $\hat{\omega}$ is the non-dimensional frequency. In our simulations, we use $\hat{\omega} = 30\pi$ and $\hat{r}_{\mathrm{core}} \approx 0.3536$, so that $\mathrm{V} \approx 2.198$. Because $V < 2.405$, the fiber is robustly single-mode.

### 4.3. Iterative Solve for the Nonlinearity.

How the nonlinear problem is solved is addressed here. It is sufficient to resort to a simple iteration scheme, where the signal and pump system is solved, and then the gain is updated, and the entire system is solved again in an iterative fashion as shown in the following algorithm: Here, $\mathfrak{u}_{l,n} = (\mathbb{E}_{l,n}, \mathbb{H}_{l,n})$ is defined to be the

---

**Algorithm 1** Simple Iterations

---

**procedure**
    $\mathfrak{u}_{l,0} = 0,\, l = s, p$
    $\Delta = \Delta_0 = 1, n = 0$
*do while* $(\Delta > \text{tol})$:
    Solve for $\mathfrak{u}_{s,n+1}$.
    Update $g_R^{p,n+1}$.
    Solve for $\mathfrak{u}_{p,n+1}$.
    Update $g_R^{s,n+1}$.
    $\Delta = \Delta_{n+1} = \frac{\|\mathfrak{u}_{s,n+1} - \mathfrak{u}_{s,n}\|}{\|\mathfrak{u}_{s,n}\|}$
*enddo*

---

electromagnetic field solutions of signal/pump ($l = s, p$) at iteration $n$ and $g_R^{l,n} = \mathbf{n}_l \Upsilon_{\mathrm{R}}^l I_k$ the corresponding gain. This process is repeated until convergence. It is worth pointing out that at each nonlinear step, a *new* (scaled) adjoint graph (test) norm is computed, which carries within it the gain contributions from the previous step:

$$\|\mathfrak{v}_n\|_{\mathfrak{Y}_n}^2 := \|\mathfrak{v}_n\|^2 + \|A_n^* \mathfrak{v}_n\|^2, \tag{4.2}$$

where

$$A_{n+1}\begin{pmatrix} \mathbb{E} \\ \mathbb{H} \end{pmatrix} := \begin{pmatrix} -(i\omega + \mathbb{P}_n) & \nabla \times \\ -\nabla \times & -i\omega \end{pmatrix} \begin{pmatrix} \mathbb{E} \\ \mathbb{H} \end{pmatrix},$$

and $\mathbb{P}_n$ is the electric polarization from the previous step. Thus, this methodology assures that the optimality properties of the ultraweak formulation are carried over at each iteration. In other words, at each step $n$, the current system of linear problems is guaranteed to be optimal. Note that by updating the test norm between each iteration, the test space is also effectively redefined between iterations. In other words, at step $n$, the test space $\mathfrak{Y}_n$ is defined by the norm 4.2, and the embedding $\mathfrak{Y}_n \hookrightarrow L^2(\Omega)$ is tacitly assumed for all $n$.

### 4.4. Optical Power Calculation.

The overall quantity of interest is the cross-sectional power through the fiber at any given $z$-value along the length of the fiber, but especially at the end of the fiber ($z = L$). Indeed, the existence of gain can be seen through the fact that energy is transferred from pump wavelength to the signal wavelength. Towards this end, one should note that the time-averaged power is computed using the (complex) Poynting vector. The (mean-squared) complex Poynting vector is defined as:

$$\mathbb{S} := \mathbb{E} \times \mathbb{H}^*,$$

where the real part $\mathbb{S}_r = \mathrm{Real}\{\mathbb{S}\}$ is the quantity of interest. Let $z = z_0$ be a position along the fiber and $\vec{n}$ be the corresponding normal vector to the cross-sectional face of the fiber at $z = z_0$. Given that most of the power in the fiber flows in the forward direction, the net power flowing in the direction determined by $\vec{n}$ through a cross-section $z = z_0$ of the fiber is computed as:

$$\mathcal{P} := \left| \int_{z=z_0} \vec{n} \cdot \mathbb{S}_r \, dS \right|.$$

In order to make rigorous mathematical sense of this term in the context of energy spaces, notice that

$$\vec{n} \cdot \mathbb{S}_r = \mathrm{Real}\{\vec{n} \cdot \mathbb{E} \times \mathbb{H}^*\} = \mathrm{Real}\{(\vec{n} \times \mathbb{E}) \cdot \mathbb{H}^*\},$$

and the last term $(\vec{n} \times \mathbb{E}) \cdot \mathbb{H}^*$ can be viewed (on the surface $z = z_0$) as a duality pairing between $H^{-1/2}(\mathrm{div}, \partial\Omega_h) \times H^{-1/2}(\mathrm{curl}, \partial\Omega_h)$. Since the UW formulation of DPG has trace variables coming from the trace spaces $H^{-1/2}(\mathrm{div}, \partial\Omega_h) \times H^{-1/2}(\mathrm{curl}, \partial\Omega_h)$, we are able to compute the power without resorting to any post-processing. Thus, the equation for $\mathcal{P}$, viewed in light of the duality pairing, has a rigorous definition.

## 5. Results

The simulation results are obtained from a single NUMA compute node with 256 GB memory and 28 cores. Initial numerical experiments on a rectangular waveguide indicated that implementing the ultraweak DPG formulation with a polynomial order $p = 5$, and with at least 4 elements per wavelength in the direction of propagation, was able to resolve the propagating wave, and so those parameters were set likewise in the fiber amplifier simulation.

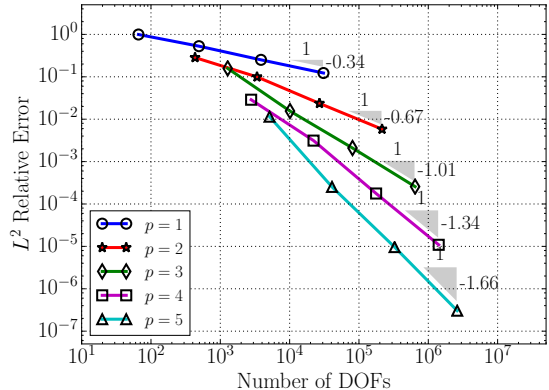

FIGURE 4. Uniform $h$-convergence rates for manufactured solution for ultra-weak DPG Maxwell formulation.

5.1. **Code Verification.** The first verification of the model consists of performing uniform $h$-convergence studies on the cylindrical core geometry comprised of curvilinear hexahedral elements. The initial mesh consists of five curvilinear hexahedra in the fiber core region. This test uses a manufactured solution of $\mathbb{E} = \sin(\omega x)\sin(\omega y)\sin(\omega z)\hat{\mathbf{e}}_x$, allowing one to find the analytical expression for the load term that is needed to produce this solution. Figure 4 depicts the expected convergence rates for the ultraweak DPG formulation (excited with $\omega = 1.001$) relative error for polynomial order $p = 1, \ldots, 5$, which theory predicts to be $-\frac{p}{3}$.

5.2. **Linear Problem.** The next verification test studies the linear case, which corresponds to setting $\mathbb{P}_l^{\text{Raman}}(\mathbb{E}_l) = 0$; in other words, this considers a simple, lossless fiber waveguide problem. In this case, only one frequency of light is needed, the signal field, denoted by $\mathbb{E}_s = \mathbb{E}$. Given that the fiber is single-mode, one ought to expect to observe the propagation of only the fundamental mode (called the $\text{LP}_{01}$ mode in scalar models) in the $x$-component of the $\mathbb{E}$ field and in the $y$-component of the $\mathbb{H}$ field. Moreover, even when the light is only seeded (introduced) into the $E_x$ component, after some distance into the fiber, one ought to expect that all of the components of the electromagnetic field acquire a non-zero value. This occurs because the light does not have to propagate perfectly in the $z$-direction, but instead is guided by total internal reflection in the core, allowing it to spread out to a small maximum angle off of the $z$-axis, which is controlled by the numerical aperture of the fiber. Only in a full vectorial model, with both electric and magnetic field components, could this phenomenon be observed.

The output $(z = L)$ images of Figs. 5-10 show that a Gaussian-shaped fundamental mode does propagate through the fiber, and that all of the field components are non-zero; though $E_x$ and $H_y$ have the largest magnitudes, as expected, since the light is seeded only into $E_x$ at $z = 0$ and there are no other polarization coupling factors that would cause energy transfer between the electromagnetic field components. Along with the real and imaginary parts of each field component (two plots on the l.h.s.), the real and imaginary parts of the cross-sectional view of the fiber parallel to the $z$-axis (two plots on the r.h.s.) are also displayed.

Another important check associated with this test is to ensure that the light does not lose power (energy) along the length of the fiber. By computing the power at various positions along the fiber (no less than one wavelength apart), the two plots of Fig. 11 are created. The start of the PML is indicated by a vertical line. These plots show that the power is (nearly) conserved as the light propagates, and, as would be expected, the solve time for the numerical model increases linearly with fiber length. Indeed, since the mesh refinements are performed anisotropically, only in the $z$-direction, the number of elements grows linearly, as does the cost of element computations. Second, the multifrontal solver has nearly linear complexity for the unidirectional anisotropically refined mesh, and the overall time grows linearly. We note that the times reported are average times for the linear solve over many runs. Specifically, fiber lengths of $L \approx 8, 16, 32, 64, 80$ are used in the test.

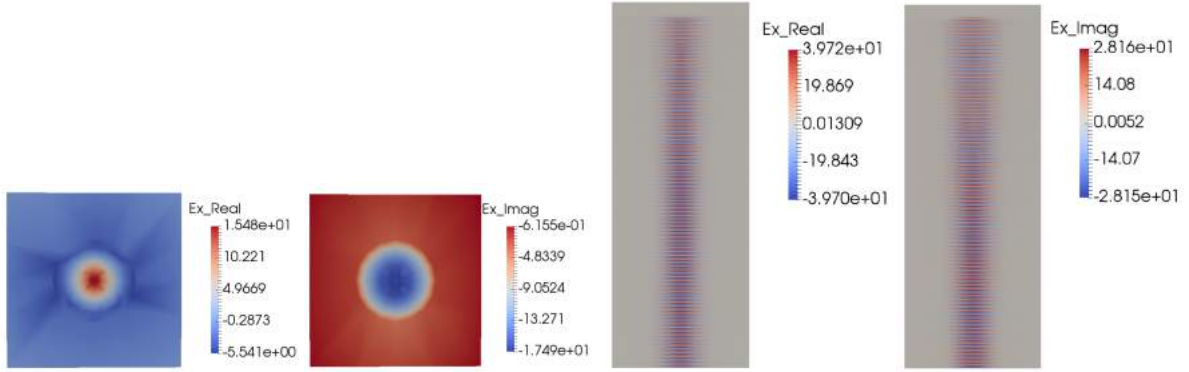

FIGURE 5. Real and Imaginary Parts of $E_x$

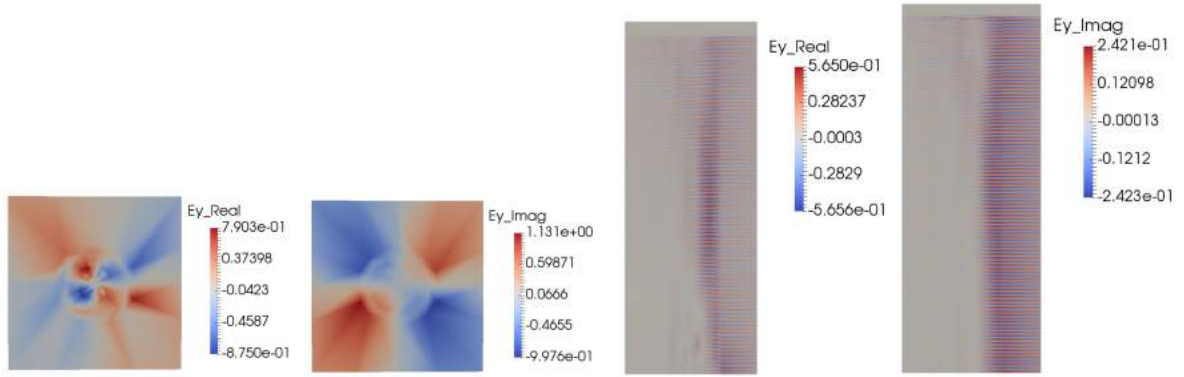

FIGURE 6. Real and Imaginary Parts of $E_y$

Since these are fibers of ultra-short lengths, elementary ray optics arguments set an upper bound in terms of the number of wavelengths required for the signal energy to settle into the physically correct solution of the waveguide. This upper bound can roughly estimated to be $\frac{r_{\text{core}}}{\tan(\text{NA})} \approx \frac{r_{\text{core}}}{(\text{NA})}$, which is ∼100-250 wavelengths for typical fibers. However, in practice, we see a physically relevant solution within $O(10)$ wavelengths, indicating that the ray optics based analysis leads to a very loose upper bound.

5.3. **Gain Problem.** The final validation of the model includes the Raman gain action along the fiber. This requires that there is both a pump wavelength and a signal wavelength, which are separated in frequency space from one another by $-13.2$ THz, corresponding to the peak Raman gain in fused silica glass. As discussed previously, the nonlinearity of the gain is handled by simple iterations. Fig. 12 illustrates the convergence of these iterations for different values of the artificial scaling $\kappa_a$. For comparison purposes, note that this test must track both frequencies of light; thus doubling the dataset size of the dependent variables. Again runs are completed for fibers of lengths: $L \approx 8, 16, 32, 64, 80$.

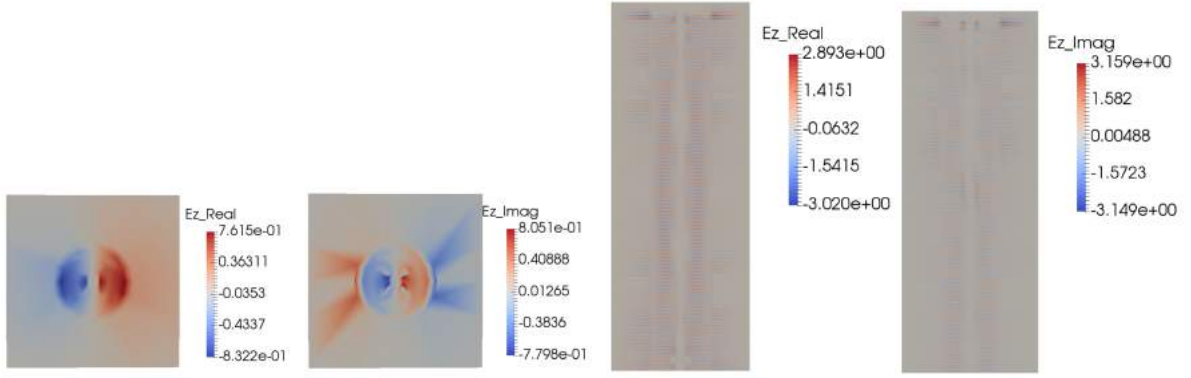

FIGURE 7. Real and Imaginary Parts of $E_z$

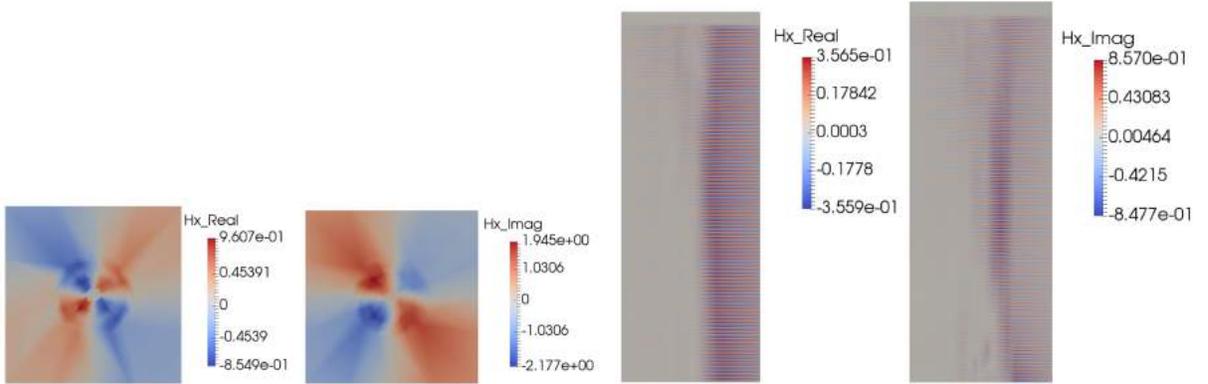

FIGURE 8. Real and Imaginary Parts of $H_x$

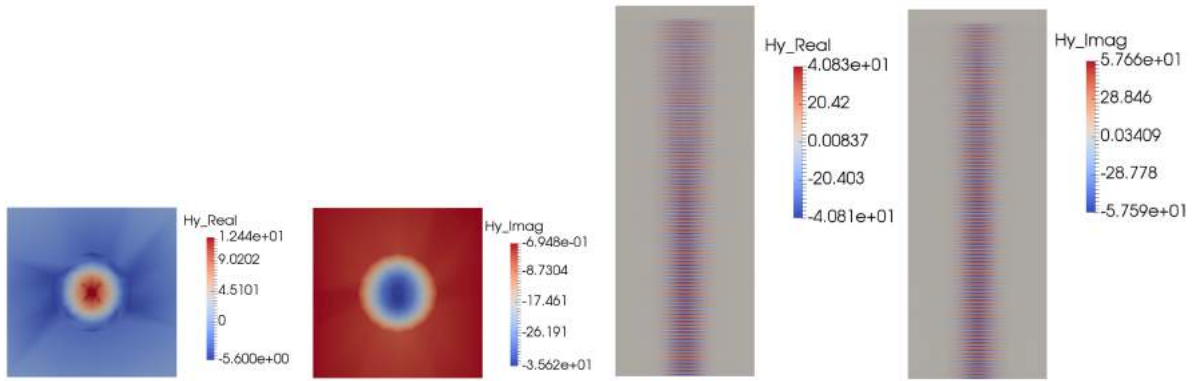

FIGURE 9. Real and Imaginary Parts of $H_y$

As with all optical nonlinearities, Raman gain per unit length increases with the intensity (irradiance) of the optical fields present in the fiber (see the coupled irradiance PDEs (2.9) or recall the gain function for Raman scattering: $g_l = \Upsilon_R^l g_R I_k$). Choosing a core-pumped amplifier, rather than a cladding-pumped amplifier configuration, ensures that the pump optical field is of a higher intensity than if it was spread out through the inner cladding and core of the fiber. The first plot of Fig. 13 depicts the pump field transferring energy to the signal field along the fiber length, as is expected from a Raman amplifier. This is for a co-pumped configuration, where both signal and pump are injected at the same fiber end ($z = 0$). The plot of Fig. 14 shows a counter-pumped configuration where the signal light is injected at $z = 0$, while the pump is introduced at $z = L$, the end of the fiber. Note that such a configuration entails the use of a separate PML for the signal and pump fields at opposite ends of the fiber. We add that such a configuration cannot be so easily modeled by a scalar BPM approach.

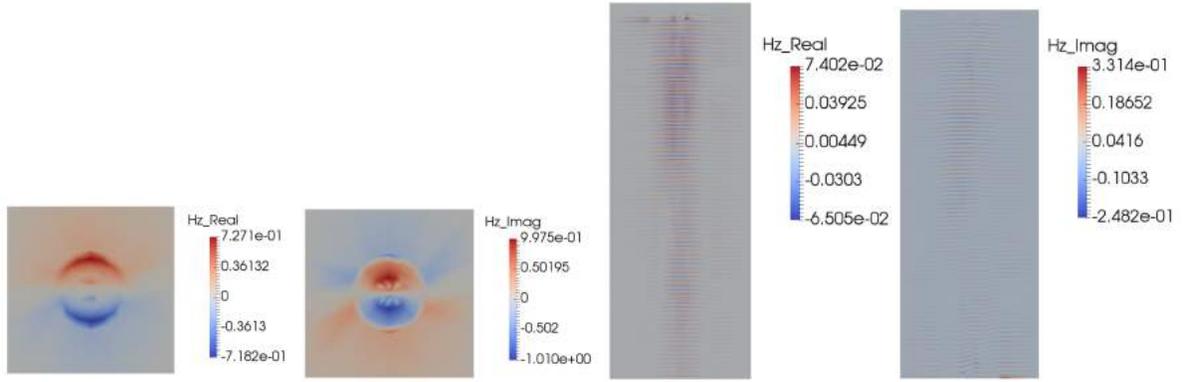

FIGURE 10. Real and Imaginary Parts of $H_z$

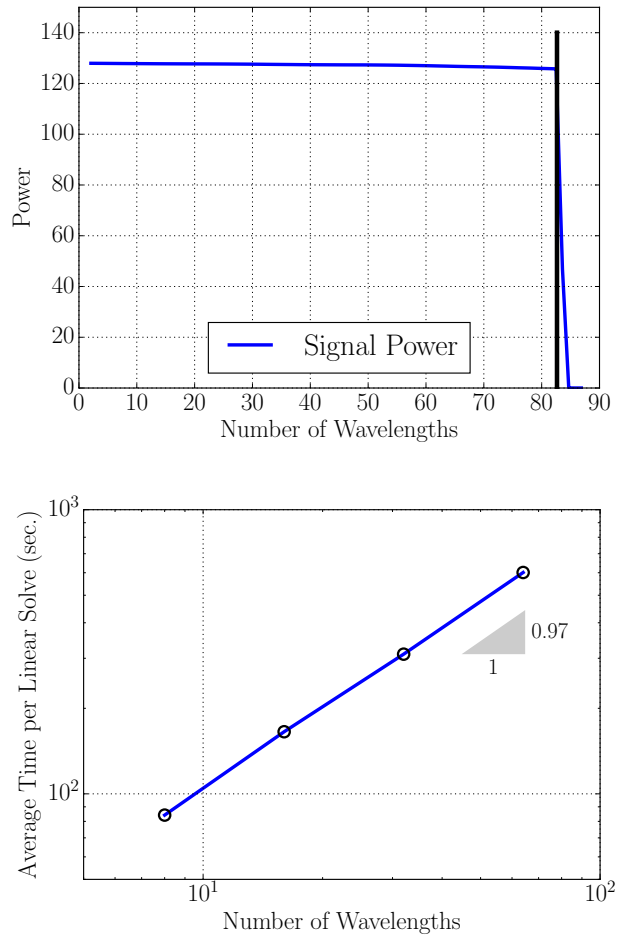

FIGURE 11. Conservation of power (top) and computational solve times for the linear model (bottom).

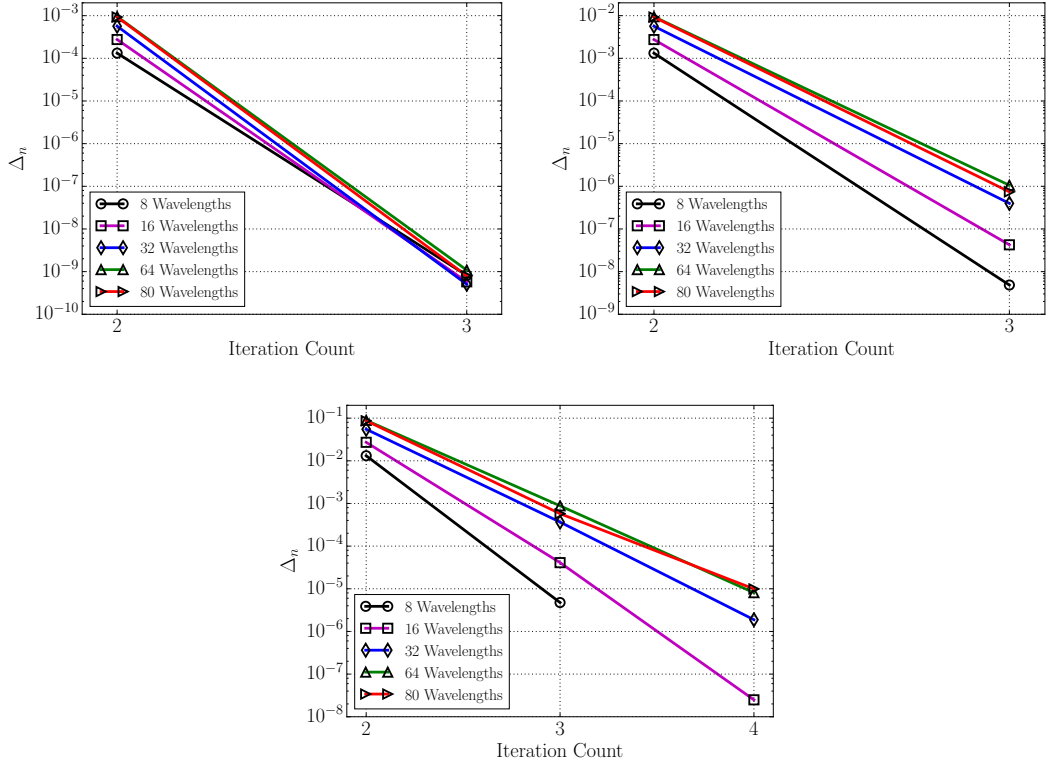

FIGURE 12. Nonlinear convergence for different $\kappa_a$.

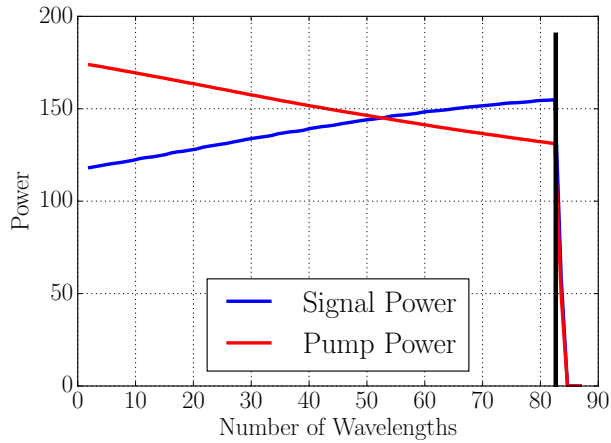

FIGURE 13. Gain for fiber of length $\approx 80$ wavelengths: co-pumped configuration with $\kappa_a = 1 \times 10^{-4}$.

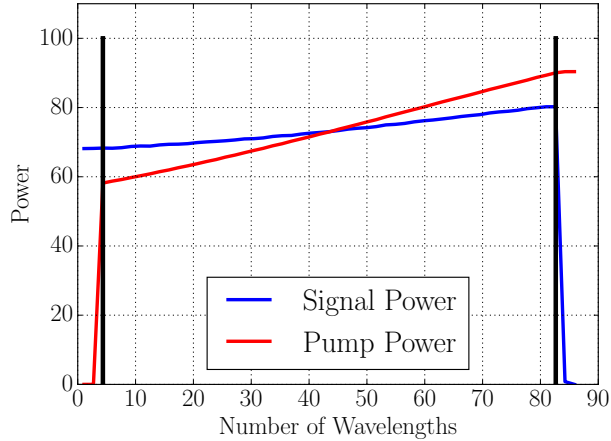

FIGURE 14. Gain for fiber of length $\approx 80$ wavelengths: counter pumped configuration with $\kappa_a = 2 \times 10^{-4}$.

The qualitative results obtained through this novel 3D vectorial DPG fiber amplifier model provide the validations needed to conclude that the methodology and implementation are sufficient for studying simple fiber amplifier configurations, and provide confidence that future efforts may prove successful in studying more complicated fiber designs under more realistic high-power operation conditions, as is the ultimate goal of this project.

## 6. Conclusions

This paper presented a unique full 3D Maxwell DPG simulation of a passive optical fiber amplifier that experiences stimulated Raman scattering. The aim was to develop computational tools for the most general model with the fewest simplifying approximations, with the intent to eventually develop a high-fidelity, multi-physics fiber model that can handle much more complex problems with realistic fiber lengths. However, in this paper, the primary interest was establishing the numerical approach and validating its feasibility by observing the qualitative characteristics of Raman gain. Towards that end, the superiority of the ultraweak DPG formulation of the coupled Maxwell system was demonstrated numerically, and was implemented in the model using a perfectly matched layer. It was successfully shown that a nonlinear iterative method was able to handle the nonlinear gain. The model verification included a convergence test, energy conservation in a linear waveguide, and qualitative gain results from a typical amplification problem. Also included was the case of a counter-pumped configuration, which is beyond the scope of most traditional scalar BPM models. This also verifies the new full vectorial electric polarization term for Raman scattering, which emphasizes the fact that the bulk Raman gain coefficient is the primary measured value available to the computer modeling team, as a practical approach to simulating Raman gain in fibers.

Several remarks are in order. First, in order to extend this model to a large-scale simulation, a significant increase in computational resources will be needed. In particular, investments must be made to develop a parallel MPI based version of the code. Second, a novel nested dissection solver, or some variant thereof, would take into advantage the possibility of static condensation of the field variables in the ultraweak formulation resulting in further optimization of the code. Third, as the complexity and/or size of the problem increases, more sophisticated nonlinear techniques may be required. From a modeling perspective, these results indicate that additional physical phenomena such as stimulated Brillouin scattering (SBS), transverse mode instability (TMI) as well as more sophisticated fiber designs and configurations, such as gain tailoring, microstructure fibers, bi-directional pumping, etc., can be accommodated within the current model. These additional modeling endeavors may require coupling DPG with other formulations [30], or implementing a coupling among various DPG formulations [29] in the Maxwell case.

# Appendices

## A. PML Details

In this appendix, we provide details of the PML implementation. Recall that the use of a PML was required due to the non-zero Raman gain term, which acts as a nonlinear conductivity. Since the model pursued in this paper is a full 3D boundary value problem (BVP) model, we must specify appropriate boundary conditions at all boundaries of the domain. We have the source located at the input ($z = 0$) end of the fiber, and PEC boundary conditions set to 0 at the external radial boundary ($r = r_{\text{cladding}}$). Given that we are modeling a fiber an arbitrary length, we must have boundary conditions of absorbing type at the (computational) exit end of the fiber ($z = z_L$). The presence of the Raman gain makes the system nonlinear, and naive absorbing impedance (or Robin) boundary conditions would induce spurious solutions. These issues can be overcome with a perfectly matched layer (PML) at the exit end of the fiber that allows the signal field to gain power and the pump field to lose power within the computational domain, while effectively setting both to 0 outside the computational domain.

Towards this end, the implementation pursued in this paper is a DPG version of a stretched coordinate PML [62], with the requirement that the stretching is done only along the $z$-direction, due to this being the direction of propagation.

### A.1. Complex stretching in $z$-direction. Let $\phi : \mathbb{R}^3 \to \mathbb{C}^3$ be a smooth invertible map with Jacobian $\mathbb{J}_{ij} = \frac{\partial \phi_i}{\partial x_j}$, where $x_i$ are the real coordinates,

$$\phi(x_1, x_2, x_3) = (\phi_1(x_1, x_2, x_3), \phi_2(x_1, x_2, x_3), \phi_3(x_1, x_2, x_3))$$

and $i, j = 1, \ldots, 3$. We let $J = |\mathbb{J}|$ denote the determinant of the Jacobian, $\mathbb{J}^{-1}$ denote its inverse and $\mathbb{J}^{-T}$ denote its inverse transpose. The map $\phi$ will be our stretching map: $\phi$ acts as identity within the computational domain, while outside, it is designed to kill outgoing waves. In our case, we have that $\phi_1, \phi_2 = 1$, since we need to stretch only the $z$-axis. Thus,

$$\mathbb{J} = \begin{pmatrix} 1 & 0 & 0 \\ 0 & 1 & 0 \\ 0 & 0 & \frac{\partial \phi_3}{\partial x_3} \end{pmatrix}.$$

The choice of the stretching function $\phi_3$ is particularly important. Given the growth of the signal field, the growth of $\phi_3$ must be commensurate so that the signal field is killed effectively.

The pump is assumed to be decaying anyway, so that it will be killed by the same PML nonetheless.

The Maxwell system can be written compactly in matrix form using the so-called Maxwell operator $A$ which acts on the electric and magnetic fields $\mathbb{E}, \mathbb{H}$. The Maxwell operator $A$ with the complex stretching included within the operator takes the form:

$$\tilde{A} \begin{pmatrix} \mathbb{E} \\ \mathbb{H} \end{pmatrix} := \begin{pmatrix} -(i\omega\epsilon + \sigma)J\mathbb{J}^{-1}\mathbb{J}^{-T} & \nabla\times \\ -\nabla\times & -i\omega\mu J\mathbb{J}^{-1}\mathbb{J}^{-T} \end{pmatrix} \begin{pmatrix} \mathbb{E} \\ \mathbb{H}, \end{pmatrix}$$

and we use the broken ultraweak formulation corresponding to $\tilde{A}$. For exponential growth, i.e., a wave of the form $e^{(a-i\omega)x_3}$, with $a > 0$, the choice of $\phi_3$ must be such that $a - \phi_3(x_1, x_2, x_3) < 0$, in order to ensure exponential decay. Figure 2, depicting the manufactured solution in the fiber waveguide is an example of a truly exponential growth, and thus $\phi_3(x_1, x_2, x_3)$ was of the form $\frac{a}{\omega}(x_3 - L)^n e^{(\frac{x_3}{L})}$ for $x_3$ values inside the PML region, where $L$ is the length of the fiber. For the simulations with Raman gain, however, such dramatic exponential growth was not observed for the fiber lengths considered, and thus $\phi_3(x_1, x_2, x_3)$ was of the form $\phi_3(x_1, x_2, x_3) = \frac{25}{\omega}\frac{(x_3-L)^3}{\beta}$, where $\beta$ is a fraction of the total length of the fiber used for the PML. For instance, $\beta = 0.2$ for the longest fiber we used.

## B. Sum Factorization Details

In this appendix, the need and efficacy of the *sum factorization* technique used in these the 3D computations is briefly reported. The exact implementation details and algorithmic break-down of this numerical integration is provided in detail in [47]. The sum factorization idea is an efficient way to dramatically reduce the time involved in the element computations (integration of element stiffness and Gram matrices) by appealing to the tensor structure of the element shapes and associated shape functions. For instance, a hexahedral element can be viewed as the tensor product of three 1D segments, and the corresponding Gaussian quadrature points of the hexahedral element can be viewed as a collection of Gaussian quadrature points of each 1D segment.

On account of the nonlinearity in this fiber amplifier problem (the Raman gain term), direct use of substructuring/templating approaches to speed-up element computations are not possible. In other words, the geometry of the fiber, which remains invariant in the longitudinal direction, cannot be used directly for recycling the stiffness or Gram matrices since they change with each nonlinear solve iteration. Other ideas such as rank-1 updates between the nonlinear iterations would also not be acceptable, since one would then be losing the true adjoint graph norm in the ultraweak formulation. Thus, the sum factorization method offers the best approach for improving the computational efficiency of the DPG implementation. Also, note that the computational improvement gained from using the sum factorization approach generally increases as the polynomial order ($p$) increases.

Sum factorization, first developed in [43] and later extended in [47], can be implemented for isoparametric elements, thus achieving the computational benefits even in problems with curvilinear geometries, which is critical for the fiber amplifier application. The ultraweak formulation and having discontinuities across elements in the $L^2$ field variables, permits the use of the sum factorization scheme for the stiffness matrices as well as the Gram matrices. In contrast, the primal formulation would have made implementing the sum factorization methodology for the stiffness matrices significantly more complicated by the need to account for the orientations of the conforming $H(\mathrm{curl})$ trial variables used in such a formulation of Maxwell's equations.

Since the hexahedral element is fully tensorial, benefits are maximized by exclusively using this element type throughout the mesh domain. While the sum factorization approach can be extended to prismatic elements, it would not only be more difficult to implement, but also one cannot guarantee the same speed-up since the prismatic element has a tensor product structure only in one direction, not all three. The use of sum factorization technique, with polynomial order $p = 5$ and hexahedral elements, in this fiber model has resulted in a computational speed-up of 80 times. It is worth noting that the many of simulations reported in this problem for longer fiber lengths, specifically over 32 wavelengths, would have been prohibitive because of the time required for element computations without this sum factorization method.

## C. DPG Theory Details

This appendix consists of two parts. First, we provide a brief, yet thorough, theoretical overview of the DPG methodology. Second, we provide numerical evidence comparing the primal and ultraweak formulations of the Maxwell system, validating that the ultraweak formulation has superior performance.

### C.1. Broken Variational Formulations and DPG.

We begin with the notion of an (abstract) broken variational formulation. A (continuous) broken variational formulation consists of a quadruple $(\mathfrak{X}, \mathfrak{Y}, \mathfrak{b}, \mathfrak{l})$, where $\mathfrak{X}, \mathfrak{Y}$ are Hilbert spaces (called the trial and test spaces respectively), $\mathfrak{b}$ is a continuous bilinear (or sesquilinear) form on $\mathfrak{X} \times \mathfrak{Y}$ and $\mathfrak{l}$ is a continuous linear (or conjugate-linear) form on $\mathfrak{Y}$. The Hilbert space $\mathfrak{X}$ is usually presented as a product of Hilbert spaces $\mathfrak{X}_0 \times \hat{\mathfrak{X}}$, while the bilinear form $\mathfrak{b}(\cdot, \cdot)$ decomposes as

$$\mathfrak{b}((\mathfrak{u}, \hat{\mathfrak{u}}), \mathfrak{v}) = \mathfrak{b}_0(\mathfrak{u}, \mathfrak{v}) + \hat{\mathfrak{b}}(\hat{\mathfrak{u}}, \mathfrak{v})$$

with $\mathfrak{b}(\cdot, \cdot), \hat{\mathfrak{b}}(\cdot, \cdot)$ being continuous bilinear (or sesquilinear) forms on $\mathfrak{X} \times \mathfrak{Y}$ and $\hat{\mathfrak{X}} \times \mathfrak{Y}$ respectively. Here, $\mathfrak{X}_0$ corresponds to the space of "field" variables while $\hat{\mathfrak{X}}$ is the interface space of trace variables. Given such a quadruple $(\mathfrak{X}, \mathfrak{Y}, \mathfrak{b}, \mathfrak{l})$ the variational problem we are interested in is the following. Find $(\mathfrak{u}, \hat{\mathfrak{u}}) \in \mathfrak{X}$ such that for all $\mathfrak{v} \in \mathfrak{Y}$, we have:

$$\mathfrak{b}((\mathfrak{u}, \hat{\mathfrak{u}}), \mathfrak{v}) = \mathfrak{l}(\mathfrak{v}). \tag{C.1}$$

The broken weak formulations of most second order equations arising in physical applications can be cast in the above abstract setting [11].

A proper understanding of the well-posedness (i.e., existence, uniqueness and stability) of such variational formulations is important to determine optimal discretization schemes. In order to determine when such an abstract broken formulation is well-posed, we make the following two assumptions:

Assumption 1. $\mathfrak{b}_0(\cdot, \cdot)$ satisfies the inf-sup condition, i.e., there exists a $\gamma > 0$ such that for all $(\mathfrak{u}, \mathfrak{v}) \in \mathfrak{X}_0 \times \mathfrak{Y}$, we have:

$$\gamma \leq \inf_{\mathfrak{u} \neq 0} \sup_{\mathfrak{v} \neq 0} \frac{|\mathfrak{b}(\mathfrak{u}, \mathfrak{v})|}{\|\mathfrak{u}\|_{\mathfrak{X}_0} \|\mathfrak{v}\|_{\mathfrak{Y}}}$$

Assumption 2. Define

$$\mathfrak{Y}_0 := \{\mathfrak{v} \in \mathfrak{Y} : \hat{\mathfrak{b}}(\hat{\mathfrak{u}}, \mathfrak{v}) = 0 \ \forall \hat{\mathfrak{u}} \in \hat{\mathfrak{X}}\}.$$

With this $\mathfrak{Y}_0$, we must ensure the triviality of the kernel $\mathfrak{Z}_0$, which is defined as

$$\mathfrak{Z}_0 = \{\mathfrak{v} \in \mathfrak{Y}_0 : \mathfrak{b}_0(\mathfrak{u}, \mathfrak{v}) = 0 \ \forall \mathfrak{u} \in \mathfrak{X}_0\}.$$

Finally, we assume $\hat{\mathfrak{b}}(\cdot,\cdot)$ satisfies the inf-sup condition, i.e., there exists a $\hat{\gamma} > 0$ such that for all $(\hat{\mathfrak{u}}, \mathfrak{v}) \in \hat{\mathfrak{X}} \times \mathfrak{Y}$, we have:

$$\hat{\gamma} \leq \inf_{\hat{\mathfrak{u}} \neq 0} \sup_{\mathfrak{v} \neq 0} \frac{|\hat{\mathfrak{b}}(\hat{\mathfrak{u}}, \mathfrak{v})|}{\|\hat{\mathfrak{u}}\|_{\hat{\mathfrak{X}}} \|\mathfrak{v}\|_{\mathfrak{Y}}}$$

Theorem 3.1 of [11] ensures that with assumptions (1) and (2), we have a well-posed variational problem corresponding to the quadruple $(\mathfrak{X}, \mathfrak{Y}, \mathfrak{b}, \mathfrak{l})$. In the sequel, we will, by abuse of notation, refer to the quadruple $(\mathfrak{X}, \mathfrak{Y}, \mathfrak{b}, \mathfrak{l})$ itself as the broken variational formulation in place of the broken variational problem defined by the quadruple. Henceforth, we assume that assumptions (1) and (2) hold and we have identified a well-posed broken variational formulation. We now discuss the DPG discretization of such a formulation.

C.1.1. *DPG Discretization.* Given a well-posed continuous broken variational formulation $(\mathfrak{X}, \mathfrak{Y}, \mathfrak{b}, \mathfrak{l})$, we come now to the task of finding an approximate (discrete) solution. Assume we have finite dimensional subspaces $\mathfrak{X}_h, \mathfrak{Y}_h$ of $\mathfrak{X}, \mathfrak{Y}$ respectively. Consider now the corresponding *discrete* variational formulation $(\mathfrak{X}_h, \mathfrak{Y}_h, \mathfrak{b}_h, \mathfrak{l}_h)$, where $\mathfrak{b}_h(\cdot,\cdot)$ and $\mathfrak{l}_h(\cdot)$ are restrictions of $\mathfrak{b}(\cdot,\cdot)$ and $\mathfrak{l}(\cdot)$ to $\mathfrak{X}_h \times \mathfrak{Y}_h$ and $\mathfrak{Y}_h$ respectively. By Babuška's theorem [3], the existence of the *discrete* inf-sup condition is *required* in order to ensure well-posedness of the *discrete variational formulation* $(\mathfrak{X}_h, \mathfrak{Y}_h, \mathfrak{b}_h, \mathfrak{l}_h)$. Regrettably, for generic finite dimensional subspaces, this cannot be guaranteed: an indiscriminate choice of discrete trial and test spaces can result in an unstable discrete formulation, *even if* the original (continuous) abstract variational formulation is well-posed. One way of defining the *ideal* DPG method is: it is the method that *achieves* the discrete inf-sup condition for *any* given discrete trial space $\mathfrak{X}_h$ by computing the so-called *optimal test space* $\mathfrak{Y}_h^{opt}$. In other words, the (ideal) DPG method answers the question, "How does one guarantee discrete stability given continuous stability?"

Since $\mathfrak{Y}$ is a Hilbert space, it comes equipped with a Riesz isometry,

$$\mathfrak{R}_{\mathfrak{Y}} : \mathfrak{Y} \to \mathfrak{Y}',$$

that maps $\mathfrak{Y}$ to its continuous dual $\mathfrak{Y}'$. Moreover, the bilinear form $\mathfrak{b}(\cdot,\cdot)$ gives rise to a linear operator $\mathfrak{B} : \mathfrak{X} \to \mathfrak{Y}'$ defined through duality between $\mathfrak{Y}$ and $\mathfrak{Y}'$ (denoted by $\langle \cdot, \cdot \rangle_{\mathfrak{Y} \times \mathfrak{Y}'}$):

$$\langle \mathfrak{B}(\mathfrak{u}), \mathfrak{v} \rangle_{\mathfrak{Y} \times \mathfrak{Y}'} := \mathfrak{b}(\mathfrak{u}, \mathfrak{v}).$$

As is well-known (see [16, 17, 11]), the trial-to-test operator $\Theta : \mathfrak{X} \to \mathfrak{Y}$, defined by $\Theta := \mathfrak{R}_{\mathfrak{Y}}^{-1} \mathfrak{B}$, when applied to the discrete trial space $\mathfrak{X}_h$, yields the optimal test space. $\Theta(\mathfrak{X}_h) = \mathfrak{Y}_h^{opt}$. With the choice of optimal test functions, one is guaranteed a unique stable solution $\mathfrak{u}_h$ of the discrete variational problem $(\mathfrak{X}_h, \mathfrak{Y}_h, \mathfrak{b}_h, \mathfrak{l}_h)$. From another viewpoint, the ideal DPG method can be viewed as a minimum residual method that minimizes the discrete operator residual i.e.,

$$\mathfrak{u}_h = \mathrm{argmin}_{\mathfrak{w}_h \in \mathfrak{X}_h} \|\mathfrak{B}(\mathfrak{w}_h) - \mathfrak{l}\|_{\mathfrak{Y}'}.$$

Indeed, upon taking the Gâteaux derivative of the operator residual and identifying the optimal test functions as $\Theta(\mathfrak{X}_h)$, we obtain the equation $\mathfrak{B}(\mathfrak{u}_h) = \mathfrak{l}$, which is equivalent to the variational problem $(\mathfrak{X}_h, \mathfrak{Y}_h, \mathfrak{b}_h, \mathfrak{l}_h)$. Finally, we can also choose to view DPG as a mixed method. We refer the reader to [14] for the details.

C.1.2. *Practical DPG.* While inversion of the test Riesz map guarantees the optimality properties of DPG, the said inversion is computationally prohibitive. Indeed, a global inversion of the test Riesz map amounts to an infinite dimensional optimization problem. A natural way out of this predicament is to consider an approximate inversion of the test Riesz map. Specifically, we fix a finite dimensional subspace $\mathfrak{Y}_r \subset \mathfrak{Y}$ with $\dim(\mathfrak{X}_h) < \dim(\mathfrak{Y}_r) < \infty$ and compute the "practical" trial-to-test operator $\Theta_r := (\mathfrak{R}_{\mathfrak{Y}}^{-1})|_{\mathfrak{Y}_r}\mathfrak{B}$, which, when applied to the discrete trial space $\mathfrak{X}_h$, yields the "practical" optimal test space. In other words, the *practical* DPG method replaces the inversion of the Riesz map $\mathfrak{R}_{\mathfrak{Y}}$ on the entire space $\mathfrak{Y}$ with inversion on the enriched space $\mathfrak{Y}_r$. A natural question at this point is, "How does this affect the optimality guarantees of DPG?" Thankfully, the answer is, not significantly. By constructing appropriate *Fortin operators* [52, 33, 11], one can control the stability of the practical DPG computations.

C.2. **Primal vs. Ultraweak Formulations.** We devote this subsection to the comparison of two (distinct) formulations of the time harmonic Maxwell system: the primal and ultraweak formulations. We assume an ansatz of the electromagnetic fields of the form $\mathbb{E}_0(x,y,z,t) = \mathbb{E}(x,y,z)e^{i\omega t}$ and $\mathbb{H}_0(x,y,z,t) = \mathbb{H}(x,y,z)e^{i\omega t}$ where $\omega > 0$ is the non-dimensionalized propagating frequency and $t$ is time. As usual, $\epsilon, \mu$ will represent non-dimensionalized electric permittivity and magnetic permeability. Moreover, we let $\sigma \in \mathbb{R}$ be a (possibly non-zero) conductivity.

C.2.1. *Primal Formulation.* The primal formulation corresponds to the case where

$$\mathfrak{X}_0 = \mathfrak{Y}_0 = H_0(\mathrm{curl}, \Omega), \hat{\mathfrak{X}} = H^{-1/2}(\mathrm{curl}, \partial\Omega_h), \mathfrak{Y} = H(\mathrm{curl}, \Omega_h).$$

The bilinear forms are:

$$\mathfrak{b}_0(\mathbb{E}, \mathbb{F}) = (\nabla \times \mathbb{E}, \nabla \times \mathbb{F})_h + ((i\omega\sigma\mu - \omega^2\epsilon\mu)\mathbb{E}, \mathbb{F})_h,$$

$$\hat{\mathfrak{b}}(\hat{\mathbb{E}}, \mathbb{F}) = \langle n \times \hat{\mathbb{E}}, \mathbb{F}\rangle_h.$$

The primal formulation is thus a broken version of the standard Bubnov-Galerkin formulation. The test space is given the standard (or mathematician's) norm, i.e.,

$$\|\mathfrak{v}\|_{\mathfrak{Y}}^2 := \|\mathbb{F}\|^2 + \|\nabla \times \mathbb{F}\|^2.$$

C.2.2. *Ultraweak Formulation.* The ultraweak formulation corresponds to the case where we integrate by parts both terms with

$$\mathfrak{X}_0 = \mathbb{L}^2(\Omega) \times \mathbb{L}^2(\Omega), \hat{\mathfrak{X}} = H^{-1/2}(\mathrm{curl}, \partial\Omega_h) \times H^{-1/2}(\mathrm{curl}, \partial\Omega_h),$$

$$\mathfrak{Y}_0 = H(\mathrm{curl}, \Omega) \times H_0(\mathrm{curl}, \Omega), \mathfrak{Y} = H(\mathrm{curl}, \Omega_h) \times H(\mathrm{curl}, \Omega_h).$$

Denote by $\mathfrak{u} = (\mathbb{E}, \mathbb{H}) \in \mathfrak{X}_0$, $\hat{\mathfrak{u}} = (\hat{\mathbb{E}}, \hat{\mathbb{H}}) \in \hat{\mathfrak{X}}$ and $\mathfrak{v} = (\mathbb{R}, \mathbb{S}) \in \mathfrak{Y}$. The bilinear forms corresponding to the ultraweak formulation are:

$$\begin{aligned}
\mathfrak{b}_0(\mathfrak{u}, \mathfrak{v}) &= (\mathbb{H}, \nabla \times \mathbb{R})_h - (i\omega\epsilon + \sigma)(\mathbb{E}, \mathbb{R})_h + (\mathbb{E}, \nabla \times \mathbb{S})_h + i\omega\mu(\mathbb{H}, \mathbb{S})_h, \\
&= (\mathfrak{u}, A^*\mathfrak{v})_h, \\
\hat{\mathfrak{b}}(\hat{\mathfrak{u}}, \mathfrak{v}) &= \langle n \times \hat{\mathbb{H}}, \mathbb{R}\rangle_h + \langle n \times \hat{\mathbb{E}}, \mathbb{S}\rangle_h,
\end{aligned} \tag{C.2}$$

equipped with the scaled adjoint graph norm:

$$\|\mathfrak{v}\|_{\mathfrak{Y}}^2 := \alpha\|\mathfrak{v}\|^2 + \|A^*\mathfrak{v}\|^2$$

with the scaling parameter $\alpha$ to make the norm localizable. It is well-known that for scalar wave propagation problems, the ultraweak formulation with the scaled adjoint graph has superior pre-asymptotic behaviour and we are guaranteed a robust estimate of the approximation error [54, 19]. While this was known to be the case for the Maxwell system as well, we provide now numerical evidence to support this claim.

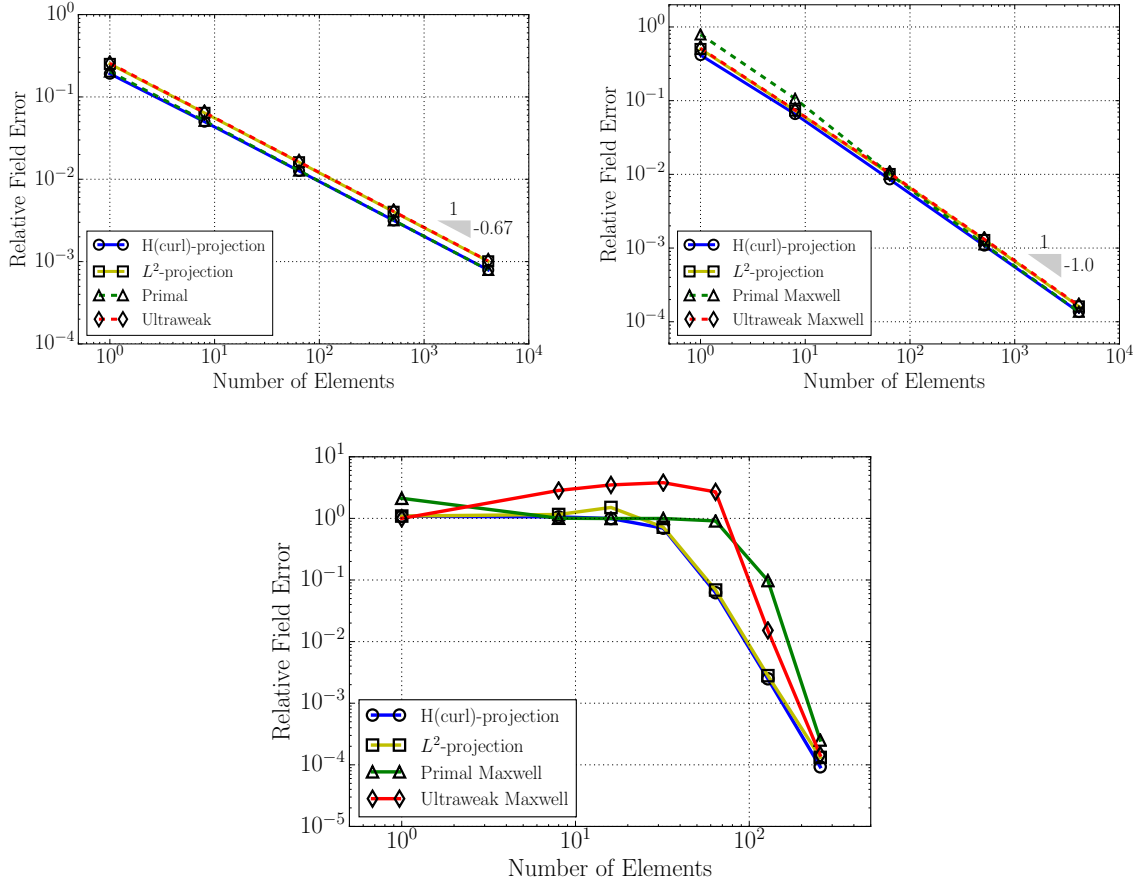

Figure 15. Energy Norm Projections and Pollution Study

C.2.3. *Energy Norm Projection and Pollution Studies.* In order to facilitate numerical comparisons of the primal and ultraweak formulations, we consider two regimes of operations. The pre- and post-asymptotic regimes. In the pre-asymptotic regime, the propagating wave is not resolved (i.e., there are not enough degrees of freedom to capture the propagation) while in the post-asymptotic regime, the wave is fully resolved. In all these cases, we assume $\sigma = 0$. Our theory indicates that the ultraweak formulation should have superior behaviour when compared with the primal formulation in both scenarios due to its use of the scaled adjoint graph norm. In the case of acoustics (Helmholtz) equation, a wavenumber explicit mathematical analysis of this behaviour is possible [19], while no such theory currently exists for the Maxwell system. The use of the scaled adjoint graph norm implies that the (ideal) unbroken ultraweak should, upon mesh refinement, deliver the $L^2$ projection in a robust fashion, while the primal has no such guarantees of convergence to the corresponding $H(\mathrm{curl})$ projection. Figure 15 shows some comparisons between the primal and ultraweak formulations. The pollution study addresses the cases of pre-asymptotic behaviour. In this study, we consider a rectangular waveguide $\Omega = [0, 1] \times [0, 1] \times [0, 16]$ and impose impedance boundary conditions on the face $z = 16$ for both formulations. The waveguide was excited by the fundamental transverse electric (TE) mode with a non-dimensionalized frequency $\omega = \sqrt{5}\pi$, which corresponds to 16 wavelengths in the $z$-direction. We now anisotropically refine the waveguide only in the $z$-direction, and study how each formulation behaves vis-a-vis the corresponding energy norm projections with polynomial order $p = 5$. We note that with $p = 5$, we required roughly 4 elements per wavelength to resolve the wave. The choice of $p = 5$ was not arbitrary, yet, lower polynomial order will require significantly more refinements to achieve the same error levels. As expected, we see that the ultraweak formulation has superior pre-asymptotic behaviour. Indeed, the error of the ultraweak formulation coincides with the $\mathbb{L}^2$ projection error earlier, while the primal formulation is farther away from the $H(\text{curl})$ projection on the same mesh. In the energy norm projection study, we study the post-asymptotic behaviour. In this case $\Omega = [0, 1]^3$ and we use polynomial orders $p = 2, 3$ and perform uniform mesh refinements. Again, we see that the ultraweak formulation "catches up" to the $\mathbb{L}^2$ projection earlier (in terms of number of refinements) than the primal with the $H(\text{curl})$ projection. This means that the ultraweak formulation (with the scaled adjoint graph norm) delivers a solution closer to the $\mathbb{L}^2(\Omega)$ projection in both the pre- and post-asymptotic regime. Thus, Fig. 15 provides numerical evidence that the ultraweak formulation is the better choice for the Maxwell system problem.

Institute for Computational Engineering and Sciences, The University of Texas at Austin, Austin, TX 78712, USA

*E-mail address*: `sriram.nagaraj@utexas.edu, Jacob.Grosek.1@us.af.mil,`
`socratis@ices.utexas.edu,`
`leszek@ices.utexas.edu, jmorapaz@ices.utexas.edu`